\documentclass[11pt]{article}
\usepackage{graphicx,amsmath,color,amsthm,amssymb}
\usepackage{psfrag,picture}
\usepackage{tikz}
\usepackage{hyperref}
\usepackage{subfig}
\usepackage{soul}
\usepackage{booktabs} 

\newcommand{\ba}{\begin{array}}
\newcommand{\ea}{\end{array}}
\newcommand{\be}{\begin{equation}}
\newcommand{\ee}{\end{equation}}
\newcommand{\ben}{\begin{equation*}}
\newcommand{\een}{\end{equation*}}
\newcommand{\bd}{\begin{displaymath}}
\newcommand{\ed}{\end{displaymath}}
\newcommand{\bi}{\begin{itemize}}
\newcommand{\ei}{\end{itemize}}
\newcommand{\bn}{\begin{enumerate}}
\newcommand{\en}{\end{enumerate}}
\newcommand{\pa}{\partial}
\newcommand{\f}{\frac}
\newcommand{\ci}{\cite}

\newcommand{\bs}{\boldsymbol}

\newtheorem{theorem}{Theorem}
\newtheorem{prob}{Problem}

\newtheorem{rem}{Remark}

\title{An Energy Based Discontinuous Galerkin Method for Coupled Elasto-Acoustic Wave Equations in Second Order Form}

\author{Daniel Appel\"{o} \thanks{Department of Applied Mathematics, Engineering Center, ECOT 225, 526 UCB, Boulder, CO 80309-0526, USA. \href{mailto:daniel.appelo@colorado.edu}{Email: daniel.appelo@colorado.edu}} ~ and Siyang Wang\thanks{Department
    of Mathematical Sciences, Chalmers University of Technology and University of Gothenburg, SE-412 96 Gothenburg,
    Sweden. \href{mailto:siyang.wang@chalmers.se}{Email: siyang.wang@chalmers.se}}}

\begin{document}
\maketitle



\begin{abstract}
We consider wave propagation in a coupled fluid-solid region, separated by a static but possibly curved interface. The wave propagation is modeled by the acoustic wave equation in terms of a velocity potential in the fluid, and the elastic wave equation for the displacement in the solid. At the fluid solid interface, we impose suitable interface conditions to couple the two equations. We use a recently developed, energy based discontinuous Galerkin method to discretize the governing equations in space. Both energy conserving and upwind numerical fluxes are derived to impose the interface conditions. The highlights of the developed scheme include provable energy stability and high order accuracy. We present numerical experiments to illustrate the accuracy property and robustness of the developed scheme. 
\end{abstract}



\section{Introduction}
Wave propagation in coupled fluid-solid domains can be found in a wide range of applications in science and industry. For example, in marine seismic surveying, when acoustic energy is released in the sea, pressure waves are created. The waves propagate from the sea through the seafloor into the earth and can be used to probe its interior. Wave propagation in the fluid can be modeled by the acoustic wave equation, and in the solid by the elastic wave equation. At the fluid-solid interface, the two equations are coupled together by physical interface conditions.

In most problems of practical interest waves travel for many temporal periods and over distances much longer than the wavelength. For such problems the leading order numerical error is typically dispersive in nature. The dispersion error is smaller with a high order method than with a low order method for the same amount of work \cite{Hagstrom2012,KreOli72} and high order methods are thus preferable. The challenge in deriving high order methods is to guarantee stability as well as accuracy. This is especially true when boundaries, or as in this paper, material interfaces are present between domains.

Here we focus exclusively on the time-dependent linear problem, where the interaction between fluid and solid is governed by the acoustic and elastic wave equations. We note that there is a very large body of research discussing the fluid-solid interaction (FSI) problem where the fluid is modeled as a nonlinear fluid and governed by the incompressible or compressible Navier-Stokes. We will not review numerical approximations of such models here but limit our discussion to high order accurate methods for the linear problem. Further, our goal is not to provide a complete review of all high order methods available but rather give a few representative references. 

For the approximation of linear waves, discontinuous Galerkin methods are energy stable, high order accurate and geometrically flexible. Recent developments include \cite{ChouShuXing2014,GSSwave} for the acoustic wave equation and \cite{ChungT119,KasDum2D06} for the elastic wave equation. Much work has been done to solve coupled acoustic-elastic problems by the discontinuous Galerkin method, for example 
in \ci{Wilcox2010} a high order discontinuous Galerkin method was developed for the coupled equations written in conservative, first order velocity-strain form. The upwind numerical fluxes of \ci{Wilcox2010} are derived by solving a Riemann problem at the interface. A generalization of the work \cite{Wilcox2010} to a formulation that avoids solving the Riemann problem can be found in \cite{ruichao}. A different discontinuous Galerkin method, taking the movement of the fluid into account and based on the first order velocity-stress formulation was proposed in \cite{Kaser2008}. The coupled viscoelastic-acoustic problem approximated by high order DG on polygonal and polyhedral meshes has been studied in \cite{antonietti2018high}. 

Hybridizable discontinuous Galerkin methods have also been used for acoustics and elastodynamics \cite{Nguyen2011,Stanglmeier2016}, and for the coupled problem \cite{Sheldon2016}. Spectral elements have been used to simulate the coupled problem in the time domain in  \cite{Komatitsch2000} and \cite{monkola2011numerical}. The latter work also shows how controllability methods can be used to solve the frequency domain problem by finding time-periodic solutions to the transient problem. 

There are also numerous finite difference methods for wave propagation problems in acoustic and elastic materials. Here we highlight some methods that satisfy a summation-by-parts principle, such that, when combined with suitable numerical boundary techniques, results in provably stable  schemes.  For example, high order finite difference methods were used to solve the acoustic wave equation in  \cite{Virta2014,Wang2016} and the elastic wave equation in \cite{Duru2014,Petersson2015}. In \cite{Petersson2018}, a finite difference method was proposed for coupling the elastic wave equation to the linearized Euler equations, modeling seismo-acoustic wave propagation. 

In \cite{Upwind2}, a general framework for energy-based discontinuous Galerkin methods for wave equations was developed. The equations are discretized directly in second order form without introducing auxiliary variables. The method admits a wide variety of mesh-independent energy conserving and upwind numerical fluxes. An application to the elastic wave equation was developed in \cite{el_dg_dath} and several test problems indicating high accuracy and robustness were presented. 

In this paper, we present an extension of the energy-based discontinuous Galerkin method in \cite{Upwind2,el_dg_dath} to coupled acoustic-elastic problems. We model wave propagation in the fluid and solid by the acoustic wave equation in terms of the velocity potential (see also \cite{Komatitsch2000}) and the elastic wave equation for displacements, respectively. Our focus in this paper is on the  treatment of fluid-solid interfaces. Stability is guaranteed by an energy analysis to the semi-discretized problem. %

In comparison with finite difference and spectral element methods discontinuous Galerkin methods are better suited to handle non-conforming refinement at interfaces but suffers from stricter time-stepping constraints  than finite difference methods. Compared to other DG methods like \ci{Wilcox2010,ruichao,Kaser2008,antonietti2018high} the method we present uses fewer variables while achieving the same or similar rates of convergence.  An advantage compared to other DG methods for wave equations in second order form is that our formulation allow for conservative and dissipative fluxes and is stable without the need for penalization.

The outline of the paper is as follows. In Section \ref{sec-eqn} we present the equations governing acoustic and elastic wave propagation. We show that the equations have an energy estimate when appropriate boundary conditions are imposed. In Section \ref{sec-dG} we outline the discontinuous Galerkin formulation for the coupled acoustic-elastic problem, and prove energy stability by deriving both energy conserving and upwind numerical fluxes. Numerical experiments including Snell's law problem, Scholte interface waves and parameter inversion verifying the effectiveness and robustness of the method are presented in Section \ref{sec-NE}. In Section \ref{sec-conclusion} we conclude and summarize. 

\section{The equations of acoustic and elastic waves}\label{sec-eqn}
In this section, we first present the equations governing wave propagation in a fluid domain $\Omega_f$ and a solid domain $\Omega_s$. The composite domain is denoted by $\Omega=\Omega_f\cup\Omega_s$, and the fluid-solid interface is $\Gamma=\Omega_f\cap\Omega_s$. We then give the physical conditions at the fluid-solid interface and show that the problem admits an energy estimate with these boundary conditions. As our focus in this paper is the fluid-solid interface, we consider both $\Omega_f$ and $\Omega_s$ to be half-planes and exclude discussions of other boundaries in the analysis. 

\subsection{The acoustic wave equation}
We consider acoustic waves propagating in an irrotational fluid in a simply connected domain $\Omega_f$. As the velocity is irrotational, the medium is a potential fluid and the wave motion can be described by the velocity potential $\psi$, which satisfies the acoustic wave equation
\be\label{a_eqn_tt}
\f{1}{c^2}\f{\pa^2\psi}{\pa t^2}= \nabla\cdot\nabla\psi,   \ \ \ \ (x_1,x_2) \in \Omega_f, \, t>0.
\ee
Here $c$ is the speed of sound. The velocity field can be obtained as the gradient of the velocity potential, i.e. $\nabla\psi$. The pressure $P=-\rho_f \pa\psi/\pa t$ is replaced by the scaled pressure $p=-P/\rho_f$, resulting in the system
\begin{alignat}{2}
& \f{\pa\psi}{\pa t}=p, & \ \   (x_1,x_2) \in \Omega_f, \, t>0, \label{a_eqn_t1}   \\
& \f{1}{c^2}\f{\pa p}{\pa t}=\nabla\cdot\nabla\psi, & \ \   (x_1,x_2) \in \Omega_f, \, t>0, \label{a_eqn_t2}
\end{alignat}
with the initial conditions
\be\label{a_eqn_ini}
\psi(x_1,x_2,0)=\psi_0(x_1,x_2),\ p(x_1,x_2,0)=p_{0}(x_1,x_2).
\ee

\subsection{The elastic wave equation} \label{sec:elasticwave}
Let $\rho = \rho(x_1,x_2)$ be the density and $u_i = u_i(x_1,x_2,t), \, i = 1,2$ be the displacements of the solid in the $x_1$ and $x_2$ direction, respectively. Assuming small deformations, the linear isotropic elastic wave equation, governing the time evolution of the displacements, is
\be\label{e_eqn_tt}
\rho \f{\pa^2 u_i}{\pa t^2}= \nabla \cdot \bs{\sigma_i},  \ \ \ \ i=1,2, (x_1,x_2) \in \Omega_s, \, t>0.
\ee
The functions $\bs{\sigma_i}:=[\sigma_{i1}, \sigma_{i2}]^T, i=1,2$ are composed of the elements of the stress tensor $\bs{\sigma}=[\bs{\sigma_1},\bs{\sigma_2}]$:
\be\label{sigma12}
\bs{\sigma_1} = \left(\begin{array}{c}
(2 \mu + \lambda) \f{\pa u_1}{\pa x_1} + \lambda \f{\pa u_2}{\pa x_2} \\
\mu (\f{\pa u_1}{\pa x_2} + \f{\pa u_2}{\pa x_1})  
\end{array} \right), \ \
\bs{\sigma_2} = \left(\begin{array}{c}
\mu (\f{\pa u_1}{\pa x_2} + \f{\pa u_2}{\pa x_1}) \\
\lambda \f{\pa u_1}{\pa x_1} + (2 \mu + \lambda) \f{\pa u_2}{\pa x_2}  
\end{array} \right),
\ee
expressed here using the Lam\'e parameters $\lambda= \lambda(x_1,x_2)$ and $\mu = \mu(x_1,x_2)$. Denoting the displacement vector $\bs{u}=[u_1,u_2]^T$, we can write \eqref{e_eqn_tt} in vector form 
\be\label{e_eqn_tt_v}
\rho\f{\pa^2\bs{u}}{\pa t^2}=\nabla\cdot\bs{\sigma}.
\ee

Introducing the velocity vector $\bs{v}=[v_1,v_2]^T$, we write \eqref{e_eqn_tt_v} as
\begin{alignat}{2}
\f{\pa \bs{u}}{\pa t}&=\bs{v},  &  \ \  (x_1,x_2) \in \Omega_s, \, t>0, \label{e_eqn_t1}\\
\rho\f{\pa \bs{v}}{\pa t}&=\nabla\cdot\bs{\sigma},  & \ \ (x_1,x_2) \in \Omega_s, \, t>0,  \label{e_eqn_t2}
\end{alignat}
with the initial conditions
\ben
\bs{u}(x_1,x_2,0) = \bs{u_{0}}(x_1,x_2), \ \ \bs{v}(x_1,x_2,0) = \bs{v_{0}}(x_1,x_2).
\een

\subsection{Interface conditions at the fluid-solid interface}
At the fluid-solid interface $\Gamma$, suitable physical conditions must be imposed so that the coupled problem is wellposed, i.e. a unique solution exists and the solution depends continuously on the data.  Let the outward pointing normal of $\Gamma$ be $\bs{n_f}$ and $\bs{n_s}$ for the fluid and solid region, respectively. Since at any point on $\Gamma$ the identity $\bs{n_f}=-\bs{n_s}$ holds,   we introduce the notation $\bs{n}=\bs{n_f}=-\bs{n_s}$. 

The fluid is attached to the solid without any gap implying that the  velocity in the normal direction is continuous,
\be\label{interface1}
\nabla\psi\cdot\bs{n}=\bs{v}\cdot\bs{n}.
\ee
In addition, the balance of forces can be stated as
\be\label{interface2}
p\bs{n}=\bs{\sigma}\cdot\bs{n}.
\ee
The operator $\cdot$ in \eqref{interface2} indicates that the tensor $\bs{\sigma}$ operates to the right on $\bs{n}$, returning a vector. Note that \eqref{interface2} is in fact two conditions. Let $\bs{m}$ be the tangential vector then we may write the two conditions as
\be\label{interface2a}
p = \bs{n} \cdot \bs{\sigma}\cdot\bs{n}, \ \  0 = \bs{m} \cdot \bs{\sigma}\cdot\bs{n}.
\ee
The first condition states balance of compression forces normal to the interface, and the second states that there is no shear force tangential to the interface. 

\subsection{Energy estimate for the coupled acoustic-elastic problem}
The natural energy estimate for the coupled acoustic-elastic problem is obtained by combining the energy estimate of the two domains. Considering first the fluid domain, we multiply \eqref{a_eqn_tt} by $\pa\psi/\pa t$ and integrate over $\Omega_f$. Assuming that the contributions from boundaries other than the interface vanish, an integration by parts yields
\ben
\int_{\Omega_f} \f{1}{c^2} \f{\pa\psi}{\pa t} \f{\pa^2\psi}{\pa t^2}\ dA = \int_{\Gamma} \f{\pa\psi}{\pa t}(\nabla\psi\cdot\bs{n})\ ds - \int_{\Omega_f}  \nabla\f{\pa\psi}{\pa t}\cdot\nabla\psi\ dA,
\een
which is equivalent to 
\ben
\f{d}{dt} E_a = \int_{\Gamma} \f{\pa\psi}{\pa t}(\nabla\psi\cdot\bs{n})\ ds,
\een
where $E_a$ is the acoustic energy 
\ben
E_a=\frac{1}{2}\int_{\Omega_f} \left[\f{1}{c^2} \left(\f{\pa\psi}{\pa t}\right)^2 + \nabla\cdot\nabla\psi\right]\ dA.
\een
The change of energy is thus determined by the contributions on the interface.

Similarly, we multiply the elastic wave equation \eqref{e_eqn_tt} by $\pa u_i/\pa t$ for $i=1,2$, integrate over $\Omega_s$ and sum over $i$. Again, the change of energy is determined by the contributions on the interface
\ben
\f{d}{dt} E_e=-\int_{\Gamma} \f{\pa \bs{u}}{\pa t} \cdot (\bs{\sigma}\cdot\bs{n})\ ds,
\een
with the elastic energy
\small
\begin{align*}
E_e=\frac{1}{2}\int_{\Omega_s} \left[\rho\left(\f{\pa u_1}{\pa t}\right)^2 + \rho\left(\f{\pa u_2}{\pa t}\right)^2\right. 
\left. +\lambda\left(\f{\pa u_1}{\pa x_1}+\f{\pa u_2}{\pa x_2} \right)^2 + \mu\left(\f{\pa u_1}{\pa x_2}+\f{\pa u_2}{\pa x_1} \right)^2 + 2\mu\left(\f{\pa u_1}{\pa x_1}\right)^2 + 2\mu\left(\f{\pa u_2}{\pa x_2}\right)^2 \right]\ dA.
\end{align*}
\normalsize
By using the interface conditions \eqref{interface1} and \eqref{interface2a} the total energy change in time is 
\begin{equation}
\begin{split}
\f{d}{dt}(E_a+E_e)& =  
\int_{\Gamma} \f{\pa\psi}{\pa t}(\nabla\psi\cdot\bs{n})
-\f{\pa \bs{u}}{\pa t} \cdot (\bs{\sigma}\cdot\bs{n})\ ds\\
&=\int_{\Gamma} p (\bs{v} \cdot \bs{n})
- [ (\bs{v}\cdot \bs{n}) \bs{n} + (\bs{v}\cdot \bs{m}) \bs{m} ] \cdot (\bs{\sigma}\cdot\bs{n})\ ds 
=
- \int_{\Gamma} (\bs{v}\cdot \bs{m}) \bs{m} \cdot (\bs{\sigma}\cdot\bs{n})\ ds =0.
\end{split}
\end{equation}
Therefore, the energy of the acoustic-elastic system is conserved. 
\begin{rem}
In the above derivation, we use the energy method to derive an energy estimate, which guarantees that the solution depends continuously on the data. For a general theory of wellposed boundary conditions for second order systems of partial differential equations, including the acoustic and elastic wave equation, see \cite{KOP_bc}.
\end{rem}

Both the acoustic wave equation and elastic wave equation can be derived by taking the variational derivative of the potential energy density as
\be\label{energy_form}
\alpha\f{\pa^2 u_i}{\pa t^2} = \sum_{k=1}^2 \f{\pa}{\pa x_k}\left(\f{\pa G}{\pa u_{i,k}} \right),
\ee
where $u_{i,k}=\pa u_i/\pa x_k$, $G$ is the potential energy density and $\alpha = 1/c^2$ for the scalar wave equation and $\alpha = \rho$ for the elastic wave equation. For the acoustic wave equation, we have $i=1$ and the potential energy density 
\be\label{ga}
G=\f{1}{2}\nabla\cdot\nabla u_1.
\ee
For the elastic wave equation, $i=1,2$ and the potential energy density is
\be\label{ge}
G=\f{\lambda}{2}\left(\f{\pa u_1}{\pa x_1}+\f{\pa u_2}{\pa x_2} \right)^2 + \f{\mu}{2}\left(\f{\pa u_1}{\pa x_2}+\f{\pa u_2}{\pa x_1} \right)^2 + \mu\left(\f{\pa u_1}{\pa x_1}\right)^2 + \mu\left(\f{\pa u_2}{\pa x_2}\right)^2.
\ee 
The discontinuous Galerkin method presented in the next section is motivated by the energy formulation of the governing equations. 

\section{A discontinuous Galerkin method for the coupled acoustic-elastic problem}\label{sec-dG}
In this section, we start by presenting a variational formulation of the coupled acoustic-elastic problem. We then propose numerical fluxes for the interface conditions \eqref{interface1} and \eqref{interface2a} that lead to a discrete energy estimate ensuring stability.  

We use the discontinuous Galerkin method developed in \cite{Upwind2} for the spatial discretization. The method is based on approximations of displacement and velocity. We therefore discretize the acoustic wave equation \eqref{a_eqn_t1}-\eqref{a_eqn_t2} and the elastic wave equation \eqref{e_eqn_t1}-\eqref{e_eqn_t2}. Crucial to the energy estimate, the energy based method uses non-standard test functions for \eqref{a_eqn_t1} and \eqref{e_eqn_t1}.

\subsection{A Galerkin variational formulation based on the energy}
Let the finite element mesh 
\[
\bar{\Omega}_f = \bigcup_i \Omega_f^i \text{\quad  and  \quad} \bar{\Omega}_s = \bigcup_j \Omega_s^j
\]
be a discretization of $\Omega_f$ and $\Omega_s$, respectively. The discretization consists of geometry-conforming and nonoverlapping quadrilaterals with piecewise smooth element boundaries. 

We use a superscript $h$ to denote the piecewise tensor product polynomial approximations of the fields. For example the velocity potential $\psi$  in the fluid, the displacement $[u_1,u_2]^T$ and velocity $[v_1,v_2]^T$ in the solid are approximated by $\psi^h$, $[u_1^h,u_2^h]^T$ and $[v_1^h,v_2^h]^T$, respectively.

Let $P_m(x)$ be a hierarchical polynomial basis in one space dimension. On a single element $\Omega_f^i$, the approximation $\psi^h$ is a tensor product polynomial in the space $\mathbb{Q}^{q_\psi}$ of degree $q_\psi$ on the reference element in the coordinate $(\xi,\eta)$. The elements of the basis are
\be\label{PiPj}
\phi_{i,j} = P_i(\xi)P_j(\eta), \ \ i,j= 0,\ldots,q_\psi.
\ee
In the derivation of variational formulations, we use subscript in the basis elements to indicate its associated variable. For example, $\phi_{\psi}$ is the basis function for $\psi$ in the form of \eqref{PiPj}. The approximations of the other variables and the associated basis functions are obtained in a similar way. 

Following \cite{Upwind2}, we test \eqref{a_eqn_t1} against $\nabla\cdot\nabla\phi_\psi$ on an element $\Omega_f^i$. After an integration by parts we find 
\begin{align*}
0 = \int_{\Omega_f^i} \left(\f{\pa\psi^h}{\pa t}-p^h\right) \nabla\cdot\nabla\phi_\psi\ dA 
=-\int_{\Omega_f^i} \nabla \left(\f{\pa\psi^h}{\pa t}-p^h\right) \cdot \nabla\phi_\psi \ dA +  \int_{\pa\Omega_f^i} \left(\f{\pa\psi^h}{\pa t}-p^h\right) (\nabla\phi_\psi\cdot\bs{n}) \ ds.
\end{align*}
By adding the penalty term 
\ben
\int_{\pa\Omega_f^i} \left(p^*-\f{\pa\psi^h}{\pa t}\right) (\nabla\phi_\psi\cdot\bs{n})\ ds
\een
to the right hand side we obtain the variational formulation for \eqref{a_eqn_t1}
\ben
\int_{\Omega_f^i} \nabla \left(\f{\pa\psi^h}{\pa t}-p^h\right) \cdot \nabla\phi_\psi\ dA= \int_{\pa\Omega_f^i} (p^*-p^h) (\nabla\phi_\psi\cdot\bs{n})\ ds.
\een
Here $p^*$ is an approximation of $\pa\psi^h/\pa t$ and $p^h$, and will be determined in the energy analysis in \mbox{Theorem \ref{theorem1}}.

We test \eqref{a_eqn_t2} against $\phi_p$,
\ben
\int_{\Omega_f^i}\left(\f{1}{c^2}\f{\pa p^h}{\pa t}-\nabla\cdot\nabla\psi^h \right)\phi_p \ dA = 0.
\een
After an integration by parts, we add a penalty term
\ben
\int_{\partial\Omega_f^i} \phi_p [(\nabla\psi\cdot\bs{n})^* - (\nabla\psi\cdot\bs{n})] \ ds,
\een
to the right-hand side, and obtain
\ben
\int_{\Omega_f^i} \left(\f{1}{c^2}\f{\pa p^h}{\pa t}\phi_p + \nabla\psi^h\cdot\nabla\phi_p\right)\ dA = \int_{\pa\Omega_f^i} \phi_p (\nabla\psi\cdot\bs{n})^*\ ds.
\een
Here $(\nabla\psi\cdot\bs{n})^*$ is an approximation of $\nabla\psi\cdot\bs{n}$ and the normal stress, and will again be determined in the energy analysis.

The variational formulations for the elastic wave equation \eqref{e_eqn_t1}-\eqref{e_eqn_t2} are obtained in a similar way. We test  the $i^{th}$ equation of \eqref{e_eqn_t1} with $\nabla\cdot\bs{\sigma_i^{\phi_u}}$, $i=1,2$,
\ben
\bs{\sigma_1^{\phi_u}} = \left(\begin{array}{c}
(2 \mu + \lambda) \f{\pa \phi_{u_1}}{\pa x_1} + \lambda \f{\pa \phi_{u_2}}{\pa x_2} \\
\mu (\f{\pa \phi_{u_1}}{\pa x_2} + \f{\pa \phi_{u_2}}{\pa x_1})  
\end{array} \right), \ \
\bs{\sigma_2^{\phi_u}} = \left(\begin{array}{c}
\mu (\f{\pa \phi_{u_1}}{\pa x_2} + \f{\pa \phi_{u_2}}{\pa x_1}) \\
\lambda \f{\pa \phi_{u_1}}{\pa x_1} + (2 \mu + \lambda) \f{\pa \phi_{u_2}}{\pa x_2}  
\end{array} \right),
\een
and obtain
\begin{align*}
0 = \int_{\Omega_s^j} \left(\f{\pa u_i^h}{\pa t}-v_i^h\right) \nabla\cdot\bs{\sigma_i^{\phi_u}} \ dA
 =  -\int_{\Omega_s^j} \nabla\left(\f{\pa u_i^h}{\pa t}-v_i^h\right) \cdot\bs{\sigma_i^{\phi_u}} \ dA - \int_{\partial\Omega_s^j}  \left(\f{\pa u_i^h}{\pa t}-v_i^h\right) (\bs{\sigma_i^{\phi_u}}\cdot\bs{n}) \ ds.
\end{align*}
Note that $\bs{n}=-\bs{n_s}$ is the source of the negative sign in the second term. The above equation can be reformulated as  
\begin{align*}
\int_{\Omega_s^j} \nabla\left(\f{\pa u_i^h}{\pa t}-v_i^h\right) \cdot\bs{\sigma_i^{\phi_u}} \ dA = - \int_{\partial\Omega_s^j}  \left(\f{\pa u_i^h}{\pa t}-v_i^h\right) (\bs{\sigma_i^{\phi_u}}\cdot\bs{n}) \ ds 
 = - \int_{\partial\Omega_s^j}  \left(v_i^*-v_i^h\right) (\bs{\sigma_i^{\phi_u}}\cdot\bs{n}) \ ds,
\end{align*}
where $v_i^*$ is an approximation of $v_i$ to be determined in the energy analysis. 

With \eqref{e_eqn_t2} tested against the standard test function $\phi_{v_1}$ and $\phi_{v_2}$ for $i=1$ and 2, respectively, we state the Galerkin variational formulation for the coupled acoustic-elastic problem.
\begin{prob}
On each element in the irrotational fluid, for all test functions 
\[
({\phi}_{\psi}, {\phi}_p) \in (\mathbb{Q}^{q_\psi} (\Omega_f^i) ) \times (\mathbb{Q}^{q_p} (\Omega_f^i)),
\]
the following variational formulation holds: 
\begin{align}
&\int_{\Omega_f^i} \nabla \left(\f{\pa\psi^h}{\pa t}-p^h\right) \cdot \nabla\phi_\psi\ dA= \int_{\pa\Omega_f^i} (p^*-p^h) (\nabla\phi_\psi\cdot\bs{n})\ ds, \label{a_eqn_t1_vf}\\
&\int_{\Omega_f^i} \left(\f{1}{c^2}\f{\pa p^h}{\pa t}\phi_p + \nabla\psi^h\cdot\nabla\phi_p\right) \ dA = \int_{\pa\Omega_f^i} \phi_p (\nabla\psi\cdot\bs{n})^* \ ds.\label{a_eqn_t2_vf}
\end{align}

On each element in the solid, for all test functions 
\[
(\bs{{\phi}_u}, \bs{{\phi}_v}) \in (\mathbb{Q}^{q_u} (\Omega_s^j) )^2 \times (\mathbb{Q}^{q_v} (\Omega_s^j))^2,
\]
the following variational formulation holds for $i=1,2$: 
\begin{align}
&\int_{\Omega_s^j} \nabla\left(\f{\pa u_i^h}{\pa t}-v_i^h\right)\cdot\bs{\sigma_i^{\phi_u}}\ dA = -\int_{\pa\Omega_s^j} (v_i^*-v_i^h)(\bs{\sigma_i^{\phi_u}}\cdot\bs{n}) \ ds, \label{e_eqn_t1_vf}\\
&\int_{\Omega_s^j} \left(\rho\f{\pa v_i^h}{\pa t}\phi_{v_i}+\bs{\sigma}_i^h\cdot\nabla\phi_{v_i}\right)\ dA = -\int_{\pa\Omega_s^j} \phi_{v_i} (\bs{\sigma_i}\cdot \bs{n})^*\ ds.\label{e_eqn_t2_vf}
\end{align}
The test function $\bs{\sigma_i^{\phi_u}}$ in \eqref{e_eqn_t1_vf} is obtained by replacing the unknown variables in \eqref{sigma12} by the associated test functions. The star variables $v_i^*$ and $(\bs{\sigma_i}\cdot \bs{n})^*$ are approximations of $v_i$ and $\bs{\sigma_i}\cdot \bs{n}$, respectively.  
\end{prob}

\subsection{Augmented equations}
The first two elements $P_0$ and $P_1$ in the polynomial basis are constant and linear. If $P_0$ is used in equation \eqref{a_eqn_t1_vf} then that equation trivially reduces to $0=0$. To obtain a  new independent equation we thus replace that equation by the moment against $P_0$, i.e. 
\be\label{a_aug}
\int_{\Omega_f^i} \left(\f{\pa\psi}{\pa t}-p\right) P_0 \ dA=0.
\ee

Similarly for $P_0$ and $i=1,2$ in equation \eqref{e_eqn_t1_vf} we find that ${0}={0}$ and replace the two missing equations by taking moments against $P_0$, i.e. 
\be\label{e_aug1}
\int_{\Omega_s^j} \left(\f{\pa u_i}{\pa t}-v_i\right) P_0 \ dA=0,\ i=1,2.
\ee

When testing \eqref{e_eqn_t1} using linear test functions \eqref{e_eqn_t1_vf} can be written as
\be\label{ss}
\int_{\Omega_s^j}\begin{bmatrix}
(2\mu+\lambda)\xi_1 & \mu\xi_2 & \mu\xi_2 & \lambda\xi_1 \\
(2\mu+\lambda)\eta_1 & \mu\eta_2 & \mu\eta_2 & \lambda\eta_1 \\
\lambda\xi_2 & \mu\xi_1 & \mu\xi_1 & (2\mu+\lambda)\xi_2 \\
\lambda\eta_2 & \mu\eta_1 & \mu\eta_1 & (2\mu+\lambda)\eta_2
\end{bmatrix}
\begin{bmatrix}
\f{\pa^2 u_1}{\pa t\pa x_1} - \f{\pa v_1}{\pa x_1}\\
\f{\pa^2 u_1}{\pa t\pa x_2} - \f{\pa v_1}{\pa x_2} \\
\f{\pa^2 u_2}{\pa t\pa x_1} - \f{\pa v_2}{\pa x_1} \\
\f{\pa^2 u_2}{\pa t\pa x_2} - \f{\pa v_2}{\pa x_2} 
\end{bmatrix}\ dA=RHS,
\ee
where $\xi_i={\pa \xi}/{\pa x_i}$, $\eta_i={\pa\eta}/{\pa x_i}$ for $i=1,2$, and $RHS$ corresponds to the right-hand side of \eqref{e_eqn_t1_vf}. The four-by-four matrix in \eqref{ss} has a zero eigenvalue, thus the equations are linearly dependent. We therefore again need to change one equation to obtain a set of linearly independent equations and a positive definite mass matrix. We opt for the strategy in Section 3.2 in \cite{el_dg_dath} by replacing the first equation in \eqref{ss} by 
\be\label{replacement}
\int_{\Omega_s^j} \left(\f{\pa^2 u_1}{\pa t\pa x_2}-\f{\pa^2 u_2}{\pa t\pa x_1}\right) - \left(\f{\pa v_1}{\pa x_2}-\f{\pa v_2}{\pa x_1}\right) \ dA= 0
\ee
if $\xi_1\eta_2=0$, or replacing the second equation in \eqref{ss} by \eqref{replacement} if $\xi_1\eta_2\neq 0$.

\subsection{Numerical fluxes for the coupled acoustic-elastic problem}
The key to couple the acoustic and elastic wave equations is to determine the stared states in \eqref{a_eqn_t1_vf}-\eqref{e_eqn_t2_vf} to impose the physical conditions \eqref{interface1} and \eqref{interface2a} at the fluid-solid interface. 

\begin{theorem}\label{theorem1}
Consider two geometry-conforming elements $\Omega_f^i\in\bar\Omega_f$ and  $\Omega_s^j\in\bar\Omega_s$ with $\Gamma_{ij}=\Omega_f^i\cap\Omega_s^j$ on the fluid-solid interface. The following numerical fluxes
\begin{align}
\bs{v^*}\cdot\bs{n} =(\nabla\psi\cdot\bs{n})^* & =\tau(\bs{v^h}\cdot\bs{n}) + (1-\tau)(\nabla\psi^h\cdot\bs{n}) - \alpha (\bs{\sigma^h}\cdot\bs{n}-p^h\bs{n})\cdot\bs{n}, \label{flux1}\\
\bs{v^*}\cdot\bs{m} &= \bs{v^h}\cdot\bs{m},\label{flux2}\\
p^*=\bs{n} \cdot (\bs{\sigma}\cdot\bs{n})^* & =\tau p^h + (1-\tau)(\bs{n} \cdot \bs{\sigma^h}\cdot\bs{n}) - \beta ((\bs{v^h}-\nabla\psi^h)\cdot\bs{n}), \label{flux3}\\
\bs{m} \cdot (\bs{\sigma}\cdot\bs{n})^*&=0, \label{flux4}
\end{align}
lead to a stable discretization for any $\tau$ if $\alpha,\beta\leq 0$. In particular, the discretization is energy-conserving if $\alpha=\beta=0$ and dissipates the energy by 
\ben
\int_{\Gamma_{ij}} \alpha (\bs{n} \cdot \bs{\sigma^h}\cdot\bs{n}-p^h)^2+\beta (\nabla\psi^h\cdot\bs{n}-\bs{v^h}\cdot\bs{n})^2\ ds,
\een
if $\alpha,\beta < 0$. 
\end{theorem}
Equation \eqref{flux1} corresponds to the continuity of velocity in the normal direction. As there are two components in $\bs{v^*}$ we need an additional condition for the velocity in the solid. In \eqref{flux2} we take the internal state of the velocity in the tangential direction. The balance of forces in the normal direction is imposed by \eqref{flux3}, and  \eqref{flux4} corresponds to zero shear force tangential to the interface. 

\begin{proof}
Replacing the test functions $\phi_\psi$ in \eqref{a_eqn_t1_vf} by $\psi^h$,  $\phi_p$ in \eqref{a_eqn_t2_vf} by $p^h$ and adding the equations we find
\ben
\f{1}{2}\f{d}{d t} \int_{\Omega_f^i} |\nabla\psi^h|^2+\f{1}{c^2}(p^h)^2\ dA = \int_{\Gamma_{ij}} (p^*-p^h)(\nabla\psi^h\cdot\bs{n})+p^h(\nabla\psi\cdot\bs{n})^*\ ds.
\een
Similarly, we replace the test functions in \eqref{e_eqn_t1_vf}-\eqref{e_eqn_t2_vf} by the corresponding numerical solutions to obtain 
\begin{align*}
\f{1}{2}\f{d}{d t}& \int_{\Omega_s^i} \lambda\left(\f{\pa u_1^h}{\pa x_1}+\f{\pa u_2^h}{\pa x_2} \right)^2 + \mu\left(\f{\pa u_1^h}{\pa x_2}+\f{\pa u_2^h}{\pa x_1} \right)^2 + 2\mu\left(\f{\pa u_1^h}{\pa x_1}\right)^2 + 2\mu\left(\f{\pa u_2^h}{\pa x_2}\right)^2\ dA  \\
= &-\int_{\Gamma_{ij}} \bs{\sigma}^h\cdot\bs{n}\cdot (\bs{v^*}-\bs{v^h}) + (\bs{\sigma}\cdot\bs{n})^*\cdot \bs{v^h}\ ds.
\end{align*}
Clearly, the discrete energy change in time is determined by the contributions on $\Gamma_{ij}$,
\be\label{I1}
I =  \int_{\Gamma_{ij}} (p^*-p^h)(\nabla\psi^h\cdot\bs{n})+p^h(\nabla\psi\cdot\bs{n})^* -\bs{\sigma^h}\cdot\bs{n}\cdot (\bs{v^*}-\bs{v^h}) - (\bs{\sigma}\cdot\bs{n})^*\cdot \bs{v^h}\ ds.
\ee
Since $\bs{v^*}$ in \eqref{I1} is not present by itself in \eqref{flux1} or \eqref{flux2}, we decompose both $\bs{v^*}$  and $\bs{v}^h$ to the normal component and tangential component as
\begin{align*}
&\bs{v^*}=(\bs{v^*}\cdot\bs{n}) \bs{n} + (\bs{v^*}\cdot\bs{m}) \bs{m}, \\
&\bs{v}^h=(\bs{v}^h\cdot\bs{n}) \bs{n} + (\bs{v}^h\cdot\bs{m}) \bs{m}.
\end{align*}
By inserting the numerical fluxes \eqref{flux1}-\eqref{flux4} to \eqref{I1} we obtain (after some algebra)
\begin{align*}
I =& \int_{\Gamma_{ij}} (\tau-1) p^h \nabla\psi^h\cdot\bs{n} + (1-\tau) \nabla\psi^h\cdot(\bs{\sigma^h}\cdot\bs{n}\cdot\bs{n})\bs{n} - \beta \nabla\psi^h\cdot ((\bs{v^h}\cdot\bs{n})\bs{n}) + \beta\nabla\psi^h\cdot(\nabla\psi^h\cdot\bs{n})\bs{n}   \\
&+\tau p^h(\bs{v^h}\cdot\bs{n}) + (1-\tau)p^h(\nabla\psi^h\cdot\bs{n}) - \alpha p^h\bs{\sigma^h}\cdot\bs{n}\cdot\bs{n}+\alpha (p^h)^2 \\
&-(\tau-1)\bs{n}\cdot\bs{\sigma}^h\cdot(\bs{v}^h\cdot\bs{n})\bs{n}-(1-\tau)\bs{n}\cdot\bs{\sigma}^h\cdot(\nabla\psi^h\cdot\bs{n})\bs{n}+\alpha\bs{n}\cdot\bs{\sigma}^h\cdot(\bs{\sigma}^h\cdot\bs{n}\cdot\bs{n}\bs{n})-\alpha p^h\bs{n}\cdot\bs{\sigma}^h\cdot\bs{n}\\
&-\tau p^h \bs{v^h}\cdot\bs{n}-(1-\tau)\bs{v}\cdot(\bs{n}\cdot\bs{\sigma}^h\cdot\bs{n})\bs{n}+\beta \bs{v^h}(\bs{v^h}\cdot\bs{n})\cdot\bs{n}-\beta\bs{v^h}\cdot(\psi^h\cdot\bs{n})\bs{n} \ ds\\
=&\int_{\Gamma_{ij}} \alpha (\bs{n}\cdot\bs{\sigma}^h\cdot\bs{n}-p^h)^2+\beta (\nabla\psi^h\cdot\bs{n}-\bs{v^h}\cdot\bs{n})^2\ ds.
\end{align*}
This proves the theorem. 
\end{proof}

The generalization to geometric non-conforming elements on the interface is straightforward by summing over contributions from all elements on the interface. 

\begin{rem}
Note that in the experiments below we use the prescription for the inter-element fluxes and enforcement of boundary conditions through numerical fluxes as described in \cite{Upwind2} and \cite{el_dg_dath} for the fluid and solid, respectively.   
\end{rem}

\subsection{Implementation of the numerical fluxes at the fluid-solid interface}
From the variational formulation \eqref{e_eqn_t1_vf}-\eqref{e_eqn_t2_vf}, we see that numerical fluxes $v_1^*$, $v_2^*$, $(\bs{\sigma_1\cdot n})^*$ and $(\bs{\sigma_2\cdot n})^*$ are needed for the implementation of the method. However, these quantities are not given directly in Theorem \ref{theorem1}. We note that \eqref{flux1}-\eqref{flux2} gives a system of two equations
\begin{equation*}
\begin{pmatrix}
n_1 & n_2 \\
-n_2 & n_1
\end{pmatrix}\begin{pmatrix}
v_1^* \\
v_2^*
\end{pmatrix}=\begin{pmatrix}
\tau(\bs{v^h}\cdot\bs{n}) + (1-\tau)(\nabla\psi^h\cdot\bs{n}) - \alpha (\bs{\sigma^h}\cdot\bs{n}-p^h\bs{n})\cdot\bs{n} \\
\bs{v^h}\cdot\bs{m} 
\end{pmatrix},
\end{equation*}
where $\bs n = [n_1, n_2]^T$ and $\bs v^* = [v_1^*, v_2^*]^T$. The solution
\begin{equation*}
\begin{pmatrix}
v_1^* \\
v_2^*
\end{pmatrix}=\begin{pmatrix}
n_1 & -n_2 \\
n_2 & n_1
\end{pmatrix}\begin{pmatrix}
\tau(\bs{v^h}\cdot\bs{n}) + (1-\tau)(\nabla\psi^h\cdot\bs{n}) - \alpha (\bs{\sigma^h}\cdot\bs{n}-p^h\bs{n})\cdot\bs{n} \\
\bs{v^h}\cdot\bs{m} 
\end{pmatrix},
\end{equation*}
can be used in the implementation of the variational formulation \eqref{e_eqn_t1_vf}. Similarly, we can solve for 
$(\bs{\sigma_1\cdot n})^*$ and $(\bs{\sigma_2\cdot n})^*$ in \eqref{flux3}-\eqref{flux4}, and obtain
\begin{equation*}
\begin{pmatrix}
(\bs{\sigma_1\cdot n})^* \\
(\bs{\sigma_2\cdot n})^*
\end{pmatrix}=\begin{pmatrix}
n_1 & -n_2 \\
n_2 & n_1
\end{pmatrix}\begin{pmatrix}
\tau p^h + (1-\tau)(\bs{n} \cdot \bs{\sigma^h}\cdot\bs{n}) - \beta ((\bs{v^h}-\nabla\psi^h)\cdot\bs{n}) \\
0
\end{pmatrix}.
\end{equation*}

\section{Numerical experiments}\label{sec-NE}
In this section, we perform numerical experiments to verify the proposed method. We start by a convergence study on a Cartesian grid for standing waves, waves governed by Snell's law, and Scholte waves. We then check the convergence rate when the method is used on curvilinear grids. We also present some more applications oriented examples.

\subsection{Standing wave problem}\label{subsection_standing}
We solve the acoustic wave equation \eqref{a_eqn_t1}-\eqref{a_eqn_t2} with wave speed $c=1$ on $\Omega_f=[0,2]^2$, and the elastic wave equation \eqref{e_eqn_t1}-\eqref{e_eqn_t2} with density $\rho=1$ and Lam\'{e} parameters $\mu=\lambda=1$ on $\Omega_s=[0,2]\times [-2,0]$. The interface between $\Omega_f$ and $\Omega_s$ is $[0,2]\times 0$. To experimentally determine rates of convergence, we use an exact solution to the acoustic wave equation 
\be\label{a_sol}
\psi = \sqrt{2} \sin(k x_1+a) \sin(k x_2+b) \sin(\sqrt{2}kt+c),
\ee
and to the elastic wave equation
\begin{align}
&u_1 = \cos(k x_1+a) \sin(k x_2+b) \cos(\sqrt{2} k t+c), \label{e_sol_1}\\
&u_2 = -\sin(k x_1+a) \cos(k x_2+b) \cos(\sqrt{2} k t+c). \label{e_sol_2}
\end{align}
Here the parameter $k\neq 0$ is used to control the wavelength. The solutions satisfy the fluid-solid  interface conditions \eqref{interface1}-\eqref{interface2}. At the boundaries, we impose Dirichlet boundary conditions. The initial and boundary data are obtained by using the exact solutions. 

We set $k=\pi$, $a=b=c=-\pi/4$ making certain that the solution at the interface is not identically zero. We use a uniform Cartesian grid with $N \times N$ square elements and side length $h=2/N$ in both $\Omega_f$ and $\Omega_s$. To evolve the solution in time we use the 8th order accurate Dormand-Prince method \cite{DORMAND198019}, with the time step $dt=0.5h/q^2$ to make sure the error in the solution is dominated by the spatial discretization. The L$_2$ errors for the four variables $\psi^h, \bs{u}^h, p^h,\bs{v}^h$ are computed at $t=2\sqrt{2}$, when the waves have propagated for two temporal periods. 

In the variational formulation \eqref{a_eqn_t1_vf}-\eqref{e_eqn_t2_vf}, numerical fluxes for both element interfaces (within either fluid or solid domain) and boundary conditions are derived in \cite{Upwind2}. In our experiments we always use upwind fluxes at the interior interfaces. For fluxes on the fluid-solid interface, we consider three cases: the upwind flux corresponding to $\tau=1/2,\ \alpha=\beta=-1$, the alternating flux $\tau=0$ or 1 with dissipation $\alpha=\beta=-1$.

The solutions $\psi^h, \bs{u}^h, p^h,\bs{v}^h$ are expanded in terms of tensor product Legendre polynomials of orders $q_\psi, q_u, q_p, q_v$, respectively. We consider two choices of the orders of approximation. First, we choose $q_\psi=q_u=q$ and $q_p=q_v=q-1$, where $q=2,3,4,5,6$. This choice is motivated by the accuracy analysis of the method applied to the acoustic wave equation in one space dimension in \cite{Upwind2}, where it is shown that the rate of convergence is optimal. By optimal, we mean that the rate of convergence is one order higher than the degree of polynomial used in the approximation. We also consider the second choice when all variables are in the same approximation space $q_\psi=q_u=q_p=q_v=q$, that is, the polynomial order of $p^h$ and $\bs{v}^h$ is increased by one. 

\begin{figure}
\begin{center}
\includegraphics[width=0.4\textwidth]{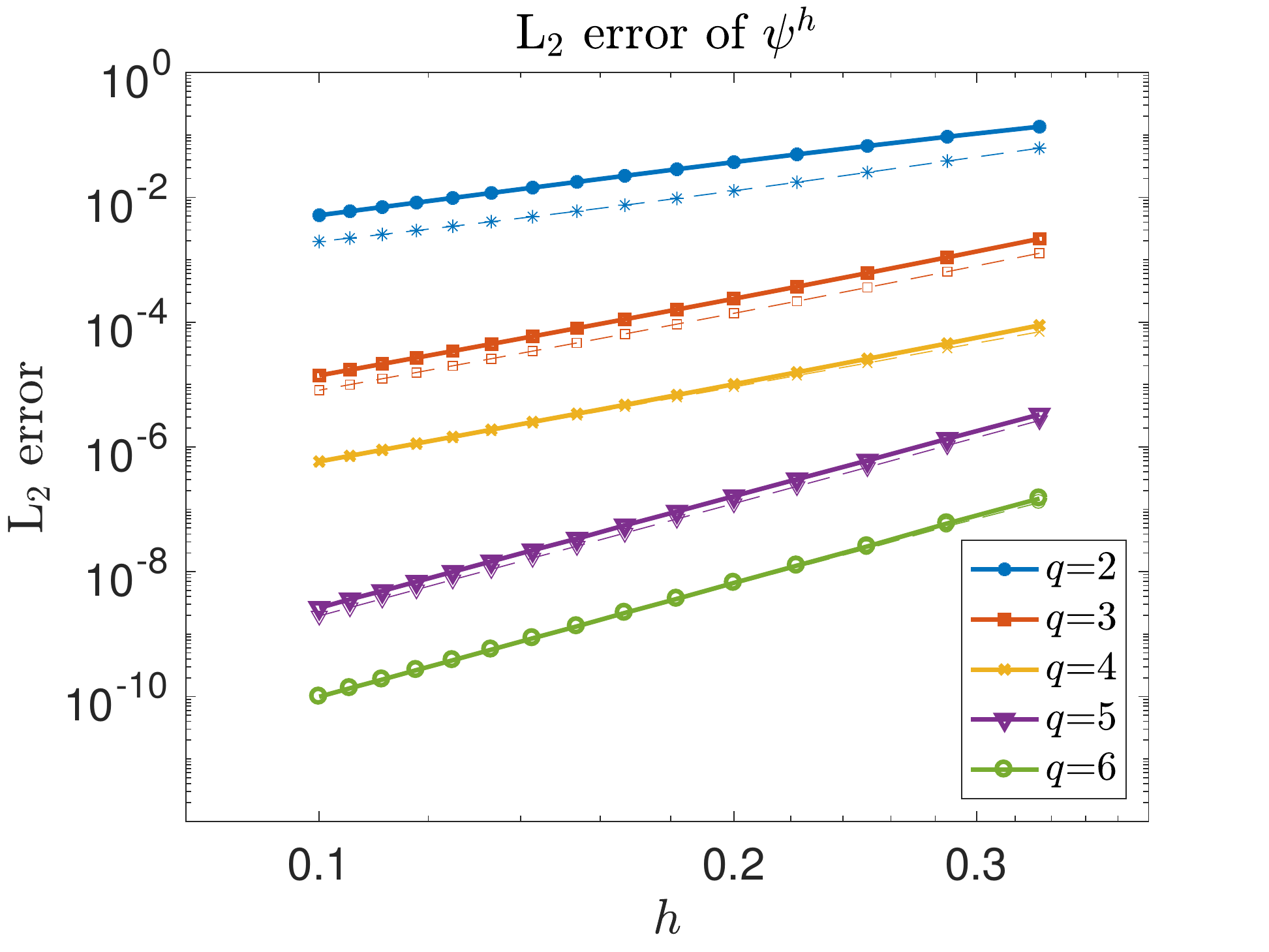}
\includegraphics[width=0.4\textwidth]{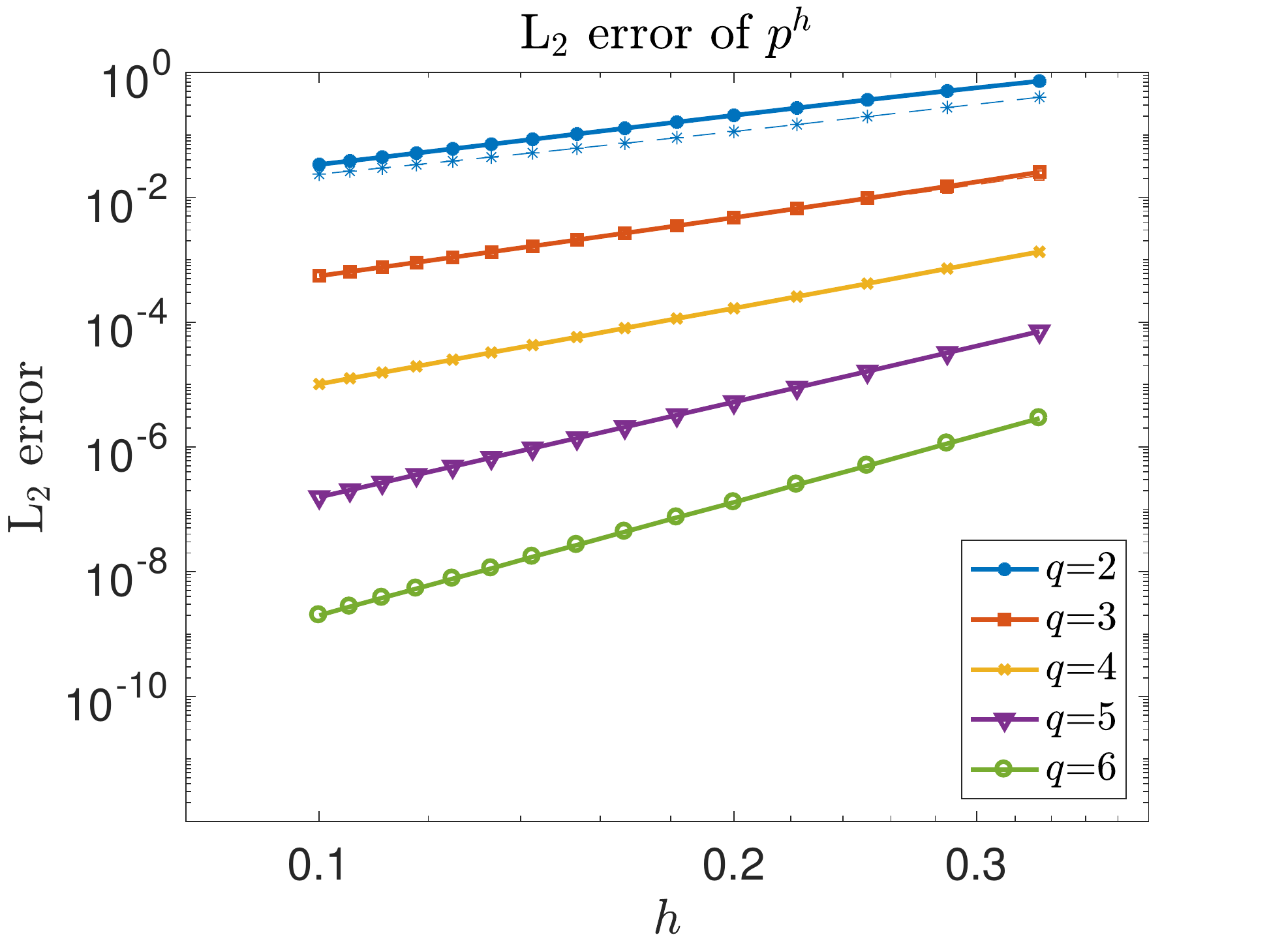}
\includegraphics[width=0.4\textwidth]{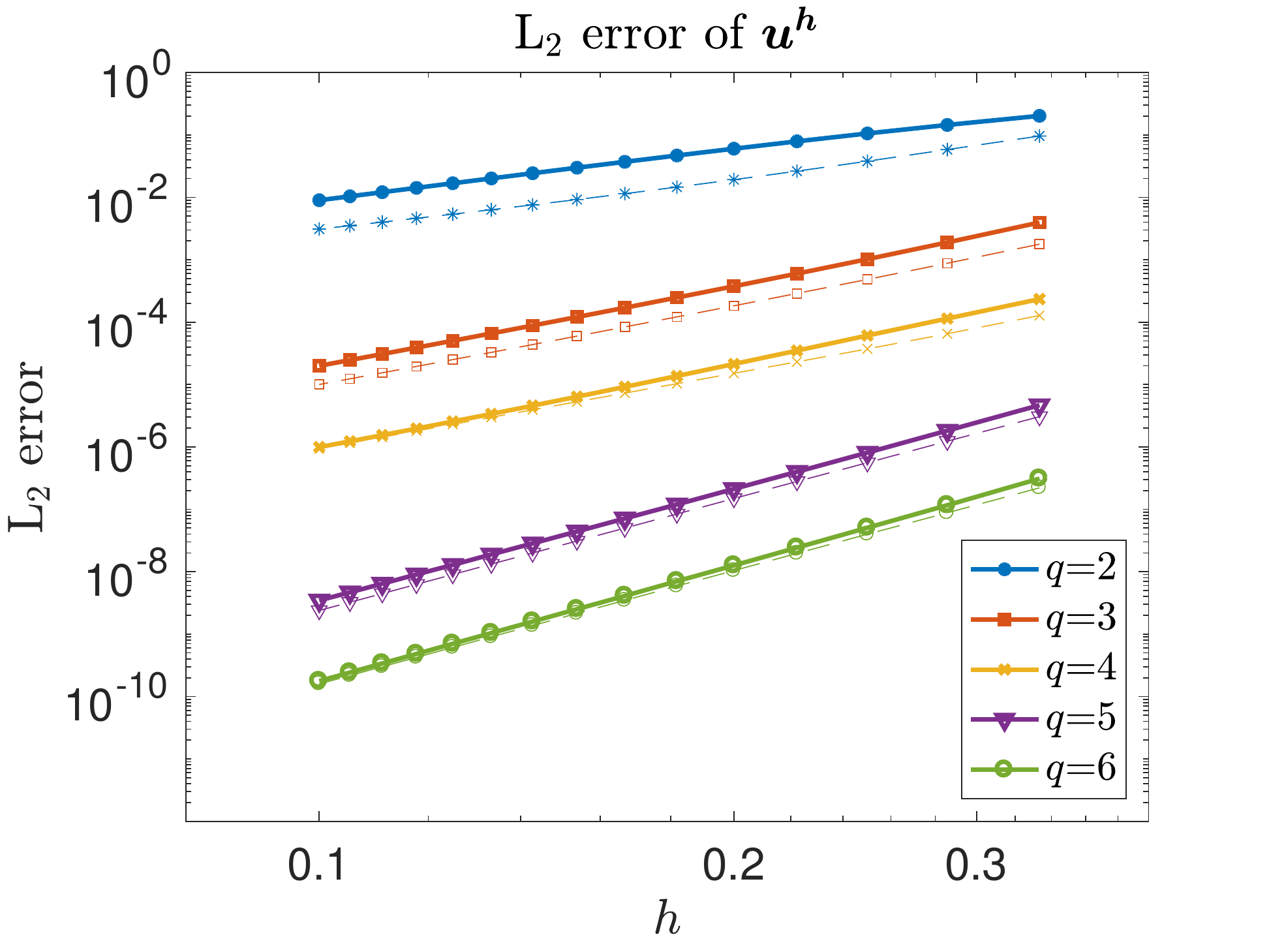}
\includegraphics[width=0.4\textwidth]{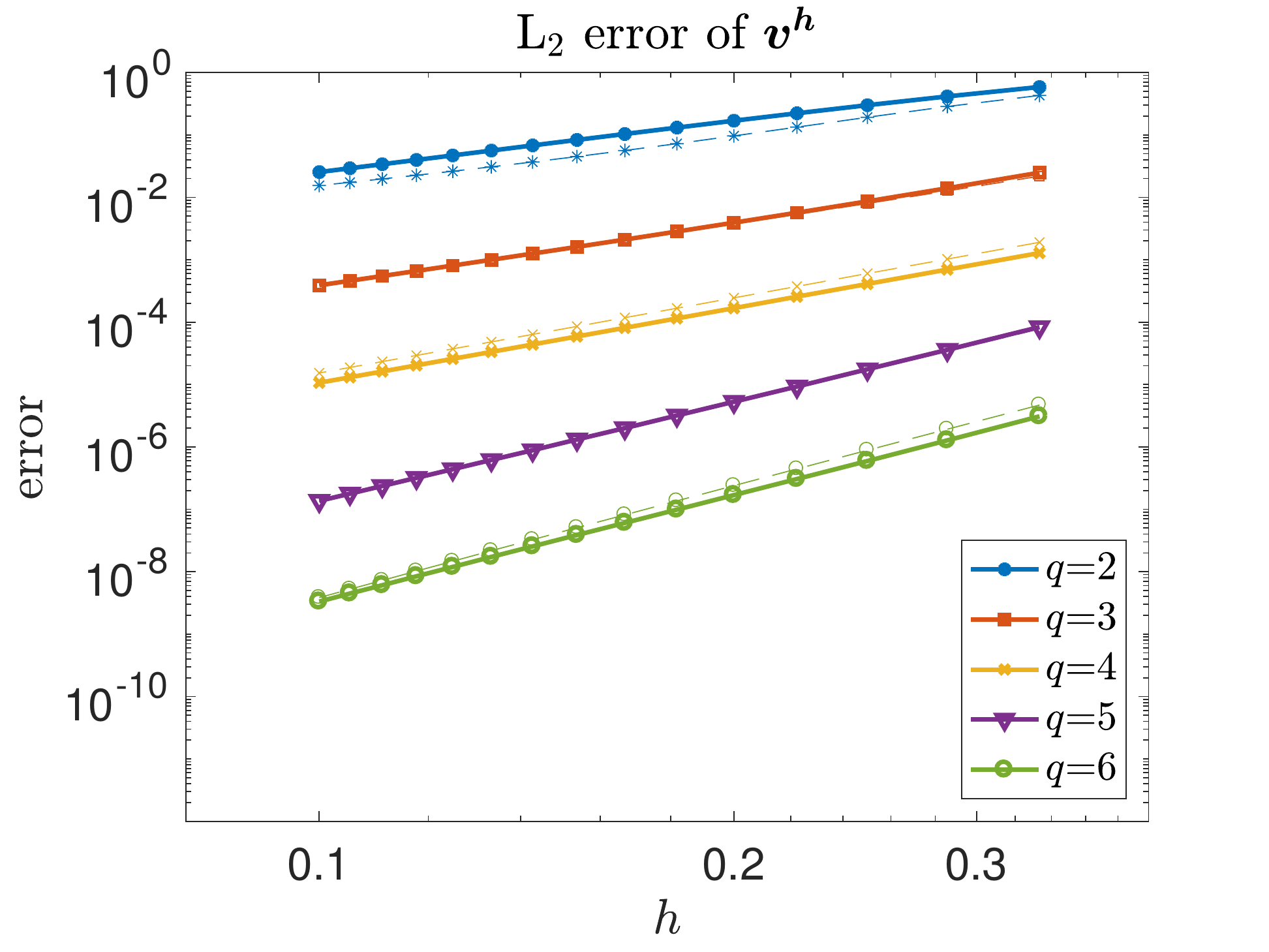}
\caption{At the top, the L$_2$ errors for the velocity potential in the fluid (left) and its time derivative (right). At the bottom, the displacement in the solid (left) and its time derivative (right). Solid lines: $q_\psi=q_u=q$ and $q_p=q_v=q-1$; dashed lines: $q_\psi=q_u=q_p=q_v=q$.}
\label{err_plot}
\end{center}
\end{figure}

To compare the above two choices of approximation spaces, we restrict to an upwind interface flux $\tau=1/2,\ \alpha=\beta=-1$. The L$_2$ error plotted in Figure \ref{err_plot} shows that the increased one order for the approximation of $p^h$ and $\bs{v}^h$ only leads to a smaller error for low order schemes with $q=2,3$, and has little influence on the error for $q>3$. In addition, rates of convergence are affected very little by this difference in approximation space, as the slopes of the two lines in the same color look almost identical.

More precisely, rates of convergence corresponding to the two choices of approximation space are shown in Table \ref{table_dq} and \ref{table_sq}, respectively. The rates of convergence are computed by the least-square fitting with the ten finest mesh refinements. In Table \ref{table_dq}, we observe that when $q$ is odd (3,5), optimal convergence is obtained for all four variables. When $q$ is even (2,4,6), convergence for $p^h$ and $\bs{v}^h$ is at least optimal. However, for the variables $\psi^h$ and $\bs{u}^h$, convergence is one order lower than optimal when $q=4$ and 6. As can be seen, the difference between the rates of convergence in Table \ref{table_dq} and \ref{table_sq} is minor. This result suggests that the first choice $q_\psi=q_u=q$, $q_p=q_v=q-1$ is better for efficiency, which is used in the following experiments in this paper.

\begin{table}
\centering
\caption{Standing wave problem: rates of convergence with $q_\psi=q_u=q$ and $q_p=q_v=q-1$, and $\tau=1/2,\ \alpha=\beta=-1$.}
\begin{tabular}{c c c c c c}
\toprule
$q$ & 2 & 3 & 4 & 5 & 6 \\
\midrule
$\psi^h$ &      2.84 & 4.07 & 4.10 & 5.94 & 6.04 \\ 
$p^h$ &         2.62 & 3.09 & 4.04 & 5.06 & 6.03 \\ 
$\bs{u}^h$ & 2.77 & 4.21 & 4.38 & 5.93 & 6.16 \\ 
$\bs{v}^h$ & 2.75 & 3.31 & 3.97 & 5.27 & 5.67 \\  \bottomrule
\end{tabular}
\label{table_dq}
\end{table} 

\begin{table}
\centering
\caption{Standing wave problem: rates of convergence with $q_\psi=q_u=q_p=q_v=q$, and $\tau=1/2,\ \alpha=\beta=-1$.}
\begin{tabular}{c c c c c c}
\toprule
$q$ & 2 & 3 & 4 & 5 & 6 \\
\midrule
$\psi^h$ &      2.66 & 4.07 & 4.01 & 5.98 & 5.99 \\ 
$p^h$ &         2.24 & 3.06 & 4.04 & 5.04 & 6.03 \\ 
$\bs{u}^h$ &  2.58 & 4.17 & 4.01 & 5.96 & 6.04 \\ 
$\bs{v}^h$ &  2.57 & 3.24 & 4.00 & 5.26 & 5.97 \\ \bottomrule
\end{tabular}
\label{table_sq}
\end{table} 

With the choice of approximation space fixed, we would like to test the influence of different numerical fluxes at the fluid-solid interface. In Table \ref{table_tau0} and \ref{table_tau1}, we show the convergence rates with alternating fluxes corresponding to $\tau=0$ and 1, respectively. We observe that the convergence rate is affected little by the choice of $\tau$.

\begin{table}
\centering
\caption{Standing wave problem: rates of convergence with $q_\psi=q_u=q$ and $q_p=q_v=q-1$, $\tau=0,\ \alpha=\beta=-1$.}
\begin{tabular}{c c c c c c}
\toprule
$q$ & 2 & 3 & 4 & 5 & 6 \\
\midrule
$\psi^h$ &      2.84 & 4.06 & 4.10 & 5.94 & 6.04 \\ 
$p^h$ &         2.62 & 3.10 & 4.04 & 5.07 & 6.03 \\ 
$\bs{u}^h$ &  2.77 & 4.21 & 4.38 & 5.93 & 6.16 \\ 
$\bs{v}^h$ &  2.75 & 3.31 & 3.98 & 5.27 & 5.70 \\ \bottomrule
\end{tabular}
\label{table_tau0}
\end{table} 

\begin{table}
\centering
\caption{Standing wave problem: rates of convergence with $q_\psi=q_u=q$ and $q_p=q_v=q-1$, $\tau=1,\ \alpha=\beta=-1$.}
\begin{tabular}{c c c c c c}
\toprule
$q$ & 2 & 3 & 4 & 5 & 6 \\
\midrule
$\psi^h$ &      2.84 & 4.07 & 4.10 & 5.95 & 6.04 \\ 
$p^h$ &         2.62 & 3.09 & 4.03 & 5.05 & 6.02 \\ 
$\bs{u}^h$ & 2.76 & 4.21 & 4.38 & 5.93 & 6.15 \\ 
$\bs{v}^h$ & 2.75 & 3.32 & 3.96 & 5.28 & 5.66 \\ \bottomrule
\end{tabular}
\label{table_tau1}
\end{table} 
In addition, we have also tested the method with Dirichlet condition at the $x_1$ boundaries, and free surface condition at the $x_2$ boundaries. The errors and rates of convergence are very close to those obtained in the above experiments with Dirichlet conditions at all boundaries, indicating the robustness of the method.

\subsection{Snell's law}
With a flat acoustic-elastic interface, an analytical solution can be derived by using Snell's law, see \cite{Wilcox2010}. When an incident pressure wave in the fluid impinges on the interface, the resulting field consists of the incident pressure wave, and also a reflected pressure wave in the fluid, transmitted pressure and shear wave in the solid. Propagation angles and wavelengths of transmitted waves are different from that of the incident wave, resulting in a more difficult test problem than the standing wave case in Section \ref{subsection_standing}.

In the fluid, the velocity potential $\psi=\psi_{ip}+\psi_{rp}$ is the sum of the incident velocity potential $\psi_{i}$ and the reflected velocity potential $\psi_{r}$, where 
\begin{align*}
&\psi_{i}=-\f{A_{i}\omega}{k_{i}}\cos(k[\sin(\alpha_{i})x_1+\cos(\alpha_{i})x_2]-\omega t),\\
&\psi_{r}=-\f{A_{r}\omega}{k_{i}}\cos(k[\sin(\alpha_{r})x_1-\cos(\alpha_{r})x_2]-\omega t). 
\end{align*}
In the solid, the displacements take the form
\begin{align*}
u_1&=A_{p} \sin({\alpha_{p}})  \cos(k_{p}[\sin(\alpha_{p})x_1+\cos(\alpha_{p})x_2]-\omega t)
-A_{s} \cos(\alpha_{s}) \cos(k_{s}[\sin(\alpha_{s})x_1+\cos(\alpha_{s})x_2]-\omega t), \\
u_2&=A_{p} \cos({\alpha_{p}})  \cos(k_{p}[\sin(\alpha_{p})x_1+\cos(\alpha_{p})x_2]-\omega t)
+A_{s} \sin(\alpha_{s})  \cos(k_{s}[\sin(\alpha_{s})x_1+\cos(\alpha_{s})x_2]-\omega t).
\end{align*}

In the ansatz, the wave numbers, wave speeds, and angular frequencies are related via
\begin{equation*}
\omega=kc=k_pc_p=k_sc_s,
\end{equation*}
and Snell's law relates  propagation angles and wave speeds as  
\begin{equation*}
\f{\sin(\alpha_{i})}{c}=\f{\sin(\alpha_{r})}{c}=\f{\sin(\alpha_{p})}{c_{p}}=\f{\sin(\alpha_{s})}{c_{s}}.
\end{equation*}

By substituting the ansatz to the interface conditions  \eqref{interface1}-\eqref{interface2}, we obtain
\begin{align*}
&A_{r}=A_{i} \f{Z_{p}(\cos(2\alpha_{s}))^2+Z_{s}(\sin(2\alpha_{s}))^2-Z}{Z_{p}(\cos(2\alpha_{s}))^2+Z_{s}(\sin(2\alpha_{s}))^2+Z},\\
&A_{p}=A_{i} \f{c}{c_{p}}\f{2Z_{p}\cos(2\alpha_{s})}{Z_{p}(\cos(2\alpha_{s}))^2+Z_{s}(\sin(2\alpha_{s}))^2+Z},\\
&A_{s}=A_{i} \f{c}{c_{s}}\f{2Z_{s}\sin(2\alpha_{s})}{Z_{p}(\cos(2\alpha_{s}))^2+Z_{s}(\sin(2\alpha_{s}))^2+Z},
\end{align*}
where
\begin{equation*}
Z=\f{c}{\cos(\alpha_{i})}, \quad Z_{p}=\f{c_{p}}{\cos(\alpha_{p})}, \quad Z_{s}=\f{c_{s}}{\cos(\alpha_{s})}.
\end{equation*} 

In the experiment, we choose a unit density in both the fluid and solid. In the fluid, the wave speed $c$ and the amplitude $A_i1$ are both chosen to be 1. In the solid, by setting the Lam\'e  parameters $\mu=4$ and $\lambda=1$, we have the pressure wave speed $c_{p}=\sqrt{(\lambda+2\mu)/\rho}=3$ and the shear wave speed $c_{s}=\sqrt{\mu/\rho}=2$.  In addition, we let the angular frequency $\omega=2\pi$, and the incident wave propagation angle $\alpha_{ip}=0.2$.

We use the same solver as in Section \ref{subsection_standing}, and compute the solution at $t=2$ when the waves have propagated for two temporal periods. Motivated by the results in Table \ref{table_dq} and \ref{table_sq}, we choose approximation for $\psi^h$ and $\bs{u}^h$ one order higher than the approximation for $p^h$ and $\bs{v}^h$, i.e.  $q_\psi=q_u:=q$, $q_p=q_v:=q-1$. In particular, we use $q=3$ and 5, as optimal or higher than optimal rates of convergence are observed in Table \ref{table_dq} with these two choices. 

\begin{figure}
\begin{center}
\includegraphics[width=0.4\textwidth]{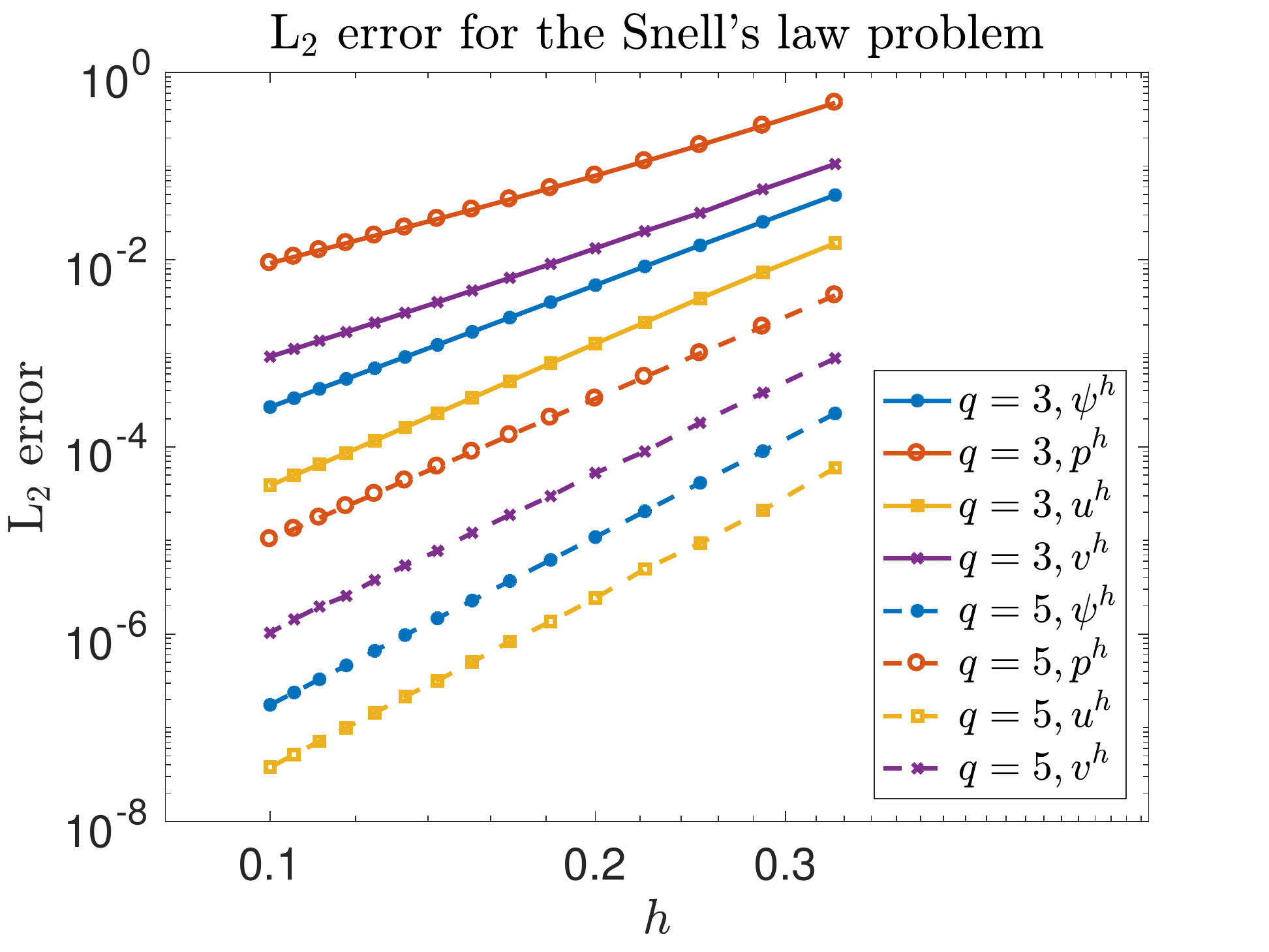}
\includegraphics[width=0.4\textwidth]{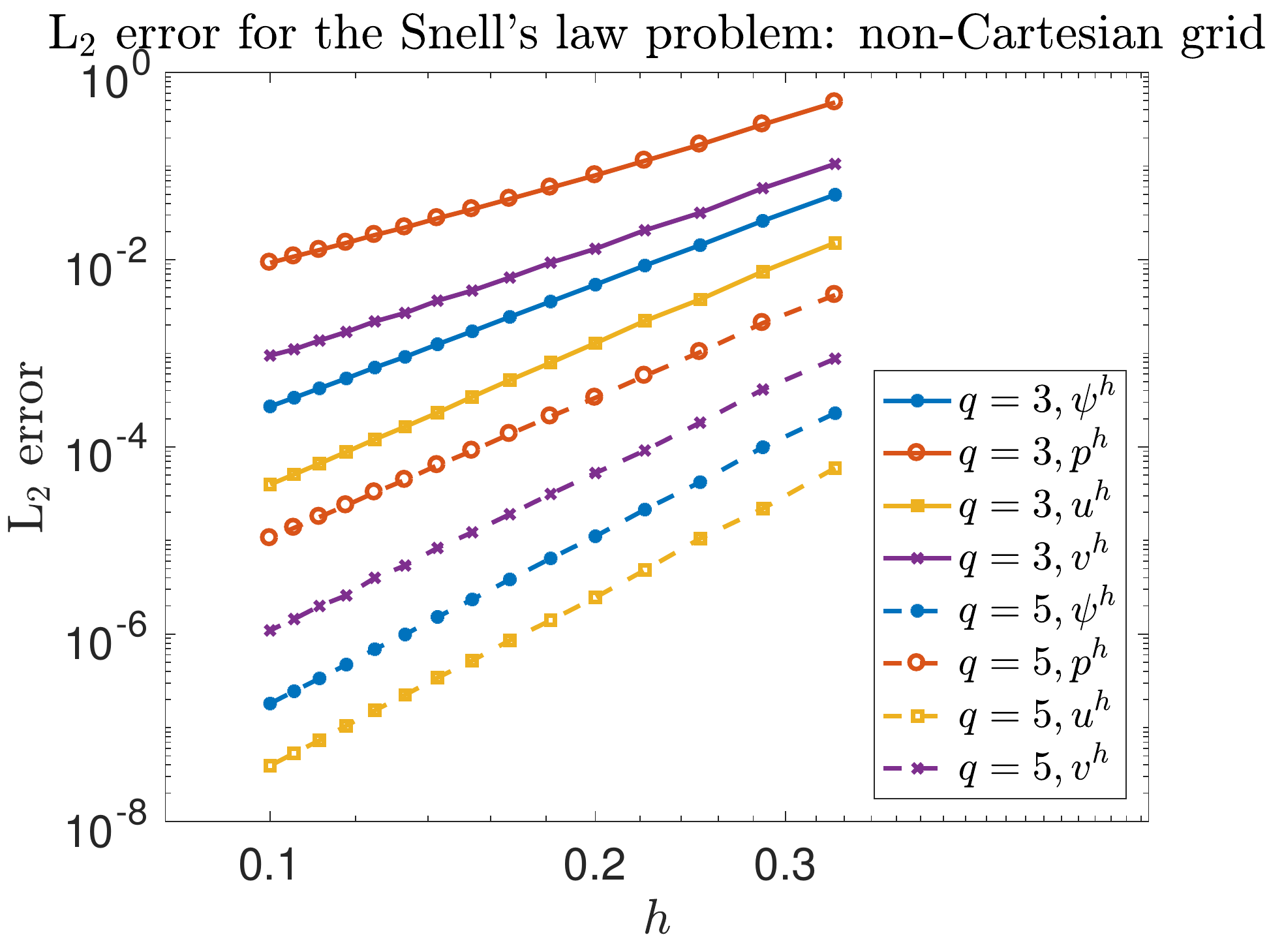}
\caption{Error plot for the Snell's law problem on a Cartesian grid (left) and on a non-Cartesian grid  (right).}
\label{err_plot_Snell}
\end{center}
\end{figure}

For the Snell's law problem, the computed rates of convergence are shown in Table \ref{table_snell}. With $q=3$, the rate of convergence is optimal for $p^h$ and higher than optimal for the other three variables. With $q=5$, we obtain an optimal rate of convergence for $\psi^h$, $p^h$, $\bs{u}^h$, and higher than optimal for $\bs{v}^h$. The error plot can be found in Figure \ref{err_plot_Snell}.

\begin{table}
\centering
\caption{Computed rates of convergence with $q_\psi=q_u=q$ and $q_p=q_v=q-1$ for the Snell's law problem on a Cartesian grid.}
\begin{tabular}{c c c c c}
\toprule
$q$ & $\psi^h$ & $p^h$ & $\bs{u}^h$  & $\bs{v}^h$  \\
\midrule
3 & 4.32 & 3.08 & 5.03 & 3.81 \\
5 & 5.96 & 5.00 & 6.01 & 5.61 \\  \bottomrule
\end{tabular}
\label{table_snell}
\end{table} 

In additional, we also test our solver for the Snell's law problem a non-Cartesian grid. The non-Cartesian grid is obtained by perturbing all interior points by $\pm 5\%$ the grid size. We plot the L$_2$ error in Figure \ref{err_plot_Snell}, and show the rate of convergence in Table \ref{table_snell_nonCartesian}. We observe that the change in the grid does not affect much the rate of convergence.
\begin{table}
\centering
\caption{Computed rates of convergence with $q_\psi=q_u=q$ and $q_p=q_v=q-1$ for the Snell's law problem on a non-Cartesian grid.}
\begin{tabular}{c c c c c}
\toprule
$q$ & $\psi^h$ & $p^h$ & $\bs{u}^h$  & $\bs{v}^h$  \\
\midrule
3 & 4.32 & 3.09 & 5.03 & 3.83 \\
5 & 5.98 & 5.02 & 6.01 & 5.62 \\  \bottomrule
\end{tabular}
\label{table_snell_nonCartesian}
\end{table}

\subsection{Snell's law between water and aluminum}
We continue with the Snell's law problem when two media have different material properties. More precisely, we consider water in the acoustic medium with a wave speed $c=1500m/s$. We use aluminum in the elastic medium with density $\rho_s = 2700 kg/m^3$, compressional wave speed $c_p=6420 m/s$, and shear wave speed $c_s=3040 m/s$. After scaling all the parameters by 1000, we solve the governing equation on a Cartesian grid in the domain $\Omega_f=[0,2]^2$ and $\Omega_s=[0,2]\times [-2,0]$. We have numerically tested that for stability the time step can be chosen 
\[
\Delta t = \frac{1.4h}{c_m (q+1.5)^2\widetilde{\rho_s}},
\]
where $c_m = \max(c,c_p)/1000=6.42$, and $\widetilde{\rho_s}=\rho_s/1000=2.7$. This amounts to a Courant number 1.4 with the scaling by the density and polynomial order. The rates of convergence shown in Table \ref{table_snell_high_contrast} are similar to the Snell's law example without high contrast parameter (see Table \ref{table_snell}), except when $q=3$  the rate for $\bs{u}^h$ drops from 5.03 to 3.83, which is slightly lower than the optimal convergence rate 4. 

\begin{figure}
\begin{center}
\includegraphics[width=0.4\textwidth]{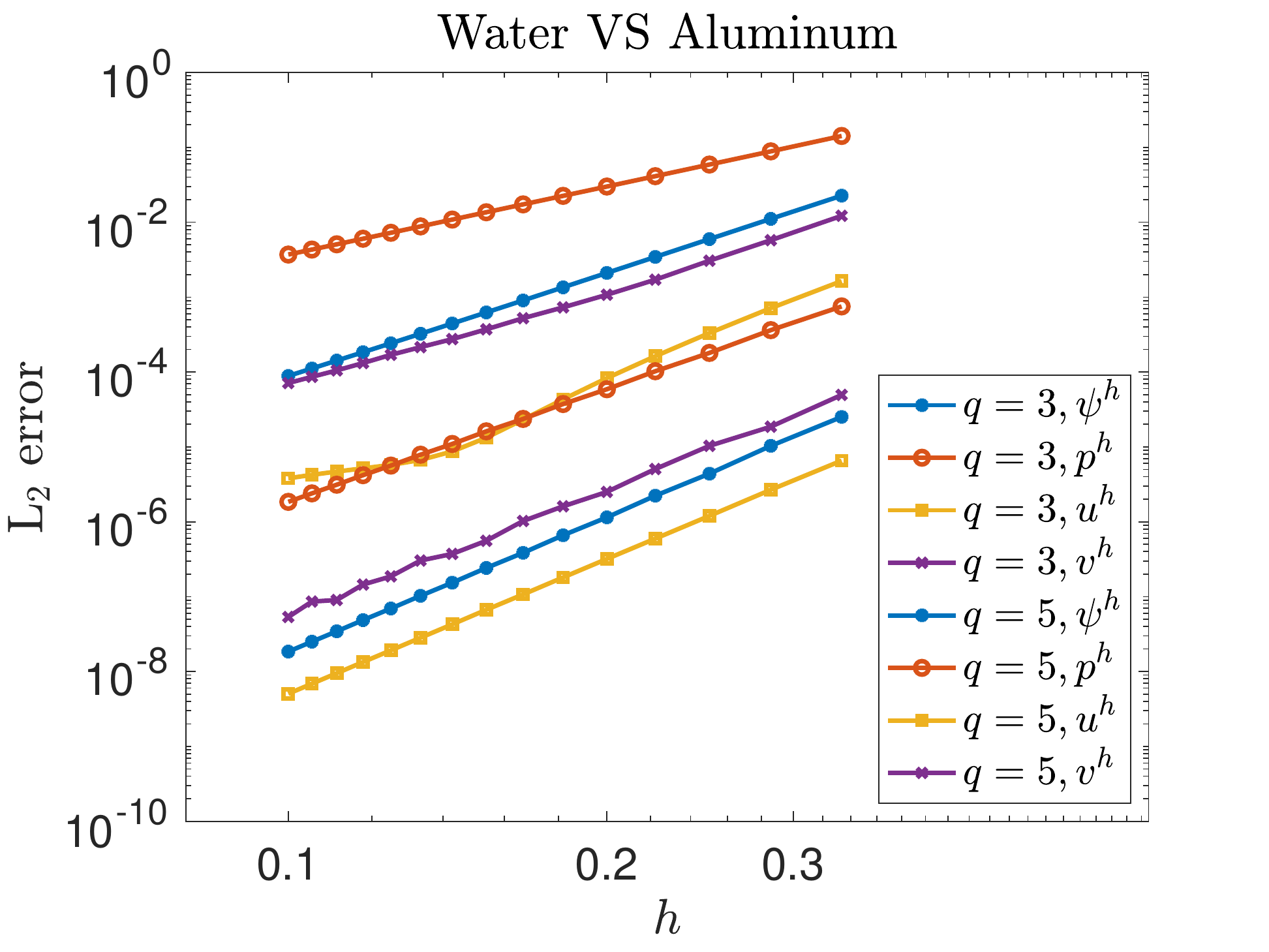}
\caption{Error plot for the Snell's law problem with water and aluminum}
\label{err_plot_Snell_WA}
\end{center}
\end{figure}

We plot the error versus the grid spacing in Figure \ref{err_plot_Snell_WA}, and observe that the L$_2$ error is at the same level of the Snell's law problem without material contrast in Figure \ref{err_plot_Snell}. 

\begin{table}
\centering
\caption{Computed rates of convergence with $q_\psi=q_u=q$ and $q_p=q_v=q-1$ for the Snell's law problem with high contrast parameters.}
\begin{tabular}{c c c c c}
\toprule
$q$ & $\psi^h$ & $p^h$ & $\bs{u}^h$  & $\bs{v}^h$  \\
\midrule
3 & 4.56 & 3.02 & 3.83 & 3.90 \\
5 & 5.97 & 5.02 & 5.98 & 5.56 \\  \bottomrule
\end{tabular}
\label{table_snell_high_contrast}
\end{table} 

\subsection{Scholte waves}
Scholte waves propagate along an acoustic-elastic interface. The waves have the highest intensity along the interface, and  decay exponentially away from the interface. This type of wave propagation is ideal for the test purpose of the proposed dG method and numerical fluxes to couple acoustic and elastic region. The velocity potential of a Scholte wave in the acoustic region $[0,2]^2$ can be written
\begin{equation*}
\psi = B_1 \omega e^{-kb_{1}x_2}\cos(kx_1-\omega t).
\end{equation*}
In the elastic region $[0,2]\times [-2,0]$, the displacements are
\begin{align*}
u_1&=(-k B_2 e^{k b_{2p}x_2}-k b_{2s}B_3e^{k b_{2s}x_2})\cos(kx_1-\omega t), \\
u_2& =(-k B_2 b_{2p} e^{k b_{2p}x_2} - k B_3 e^{k b_{2s}x_2}) \sin(kx_1-\omega t).
\end{align*}

We take material parameters $c=1$ in the acoustic wave equation, and $\rho=\lambda=\mu=1$ in the elastic wave equation. The decay rates are
\begin{equation*}
b_1 = \sqrt{1-c_s^2/c_{1}^2}, \quad b_{2p}=\sqrt{1-c_s^2/c_{2p}^2}, \quad b_{2s}=\sqrt{1-c_s^2/c_{2s}^2},
\end{equation*}
where the wave speeds in the acoustic and elastic regions
\begin{equation*}
c_1=c=1, \quad c_{2p}=\sqrt{\f{\lambda+2\mu}{\rho}}, \quad c_{2s}=\sqrt{\f{\mu}{\rho}}.
\end{equation*}
The speed of Scholte wave $c_s$ is determined by the interface conditions \eqref{interface1}-\eqref{interface2}, and the wave number is $k=\omega/c_s$. In the experiment, we use the same parameters as in \cite{Wilcox2010} and choose $c_s=0.7110017230197$, $B_1 = -0.3594499773037$, $B_2 = -0.8194642725978$, $B_3 =1$, and $\omega=2\pi$.

\begin{figure}
\begin{center}
\includegraphics[width=0.4\textwidth]{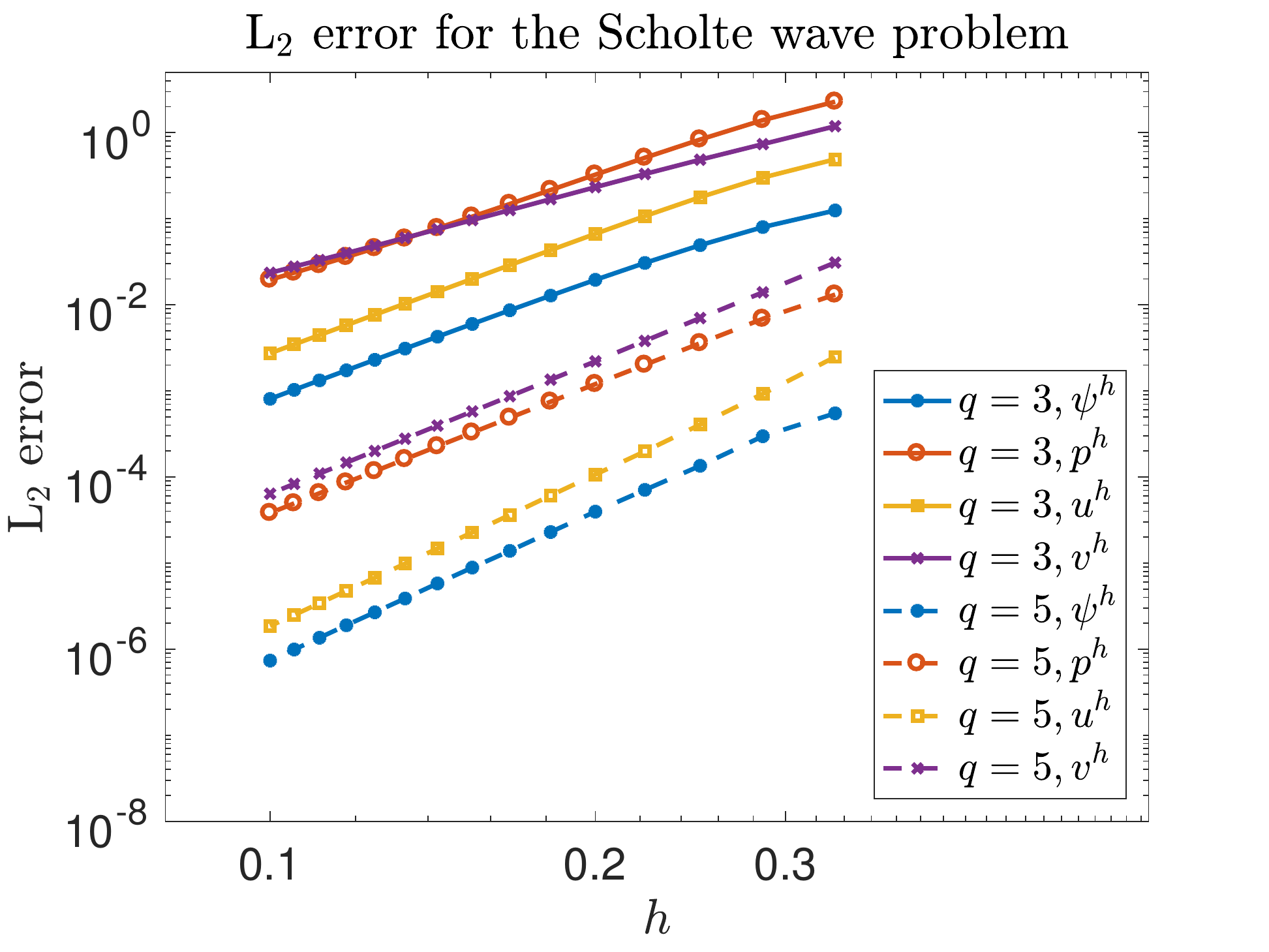}
\caption{Error plot for the Scholte wave problem.}
\label{err_plot_Scholte}
\end{center}
\end{figure}

The computed rates of convergence at $t=2$ are shown in Table \ref{table_scholte}. For $q=3$, we observe higher than optimal rates for all the four variables. For $q=5$, the rates are optimal for $p^h$ and $\bs{v}^h$, and slightly lower than optimal for $\psi^h$ and $\bs{u}^h$. We also plot the L$_2$ error versus the grid spacing in Figure \ref{err_plot_Scholte}.

\begin{table}
\centering
\caption{Computed rates of convergence with $q_\psi=q_u=q$ and $q_p=q_v=q-1$ for the Scholte wave.}
\begin{tabular}{c c c c c}
\toprule
$q$ & $\psi^h$ & $p^h$ & $\bs{u}^h$  & $\bs{v}^h$  \\
\midrule
3 & 4.63 & 4.01 & 4.61 & 3.30 \\
5 & 5.75 & 4.98 & 5.83 & 5.10 \\  \bottomrule
\end{tabular}
\label{table_scholte}
\end{table}

\subsection{Modes of two coupled annulus}
In this experiment, we apply our method on a curvilinear mesh. Consider a domain consisting of an annular fluid region $r_1 \le r \equiv \sqrt{x_1^2+x_2^2} \le r_2$ and a solid region confined in $r_0 \le r  \le r_1$. For simplicity we assume that the speed of sound in the fluid is one and that the P-velocity in the solid is also one. Precisely we take the density to be one and $\lambda = 0.5$ and $\mu = 0.25$. 

Denote by $u_r$ and $v_r$ the radial component of the solid displacement and velocity, then a solution to the elastic wave equation is 
\be
u_r(r,t) = J_1(r) \cos (t) ,\ \ v_r(r,t) = -J_1(r) \sin (t).  
\ee
Similarly a solution to the wave equation in the fluid is 
\be
\psi(r,t) = J_0(r) \sin (t), \ \ p(r,t) = J_0(r) \cos (t). 
\ee  
Above $J_n(r)$ is the Bessel function of the first kind of order $n$.

At $r = r_1$ the interface conditions become
\begin{align}
p &= (2\mu + \lambda) \frac{\partial u_r}{\partial r} + \frac{\lambda}{r_1} , \\
\frac{\partial \psi}{\partial r} &= v_r, 
\end{align} 
which reduces to the solvability condition 
\be
(2\mu + \lambda)  \frac{\partial J_1(r_1)}{\partial r}   \frac{\partial J_0(r_1)}{\partial r}
 + \frac{\lambda}{r_1} J_1(r_1) \frac{\partial J_0(r_1)}{\partial r} + \frac{\lambda}{r_1} J_1(r_1)  J_0(r_1) = 0. \label{ann_solv}
\ee
Noting that $J'_0(r) = -J_1(r)$ it is clear that $r_1$ corresponding to roots of $J_1(r)$ are also solutions to (\ref{ann_solv}). We thus set $r_1 = 7.01558666981561$ and also select  $r_0 = 3.83170597020751$ and $r_2 = 8.65372791291101$ to be roots of $J_1$ and $J_0$, respectively, so that the boundary conditions become homogenous and of Dirichlet type on the inner and outer boundaries.

\begin{figure}
\begin{center}
\includegraphics[width=0.4\textwidth]{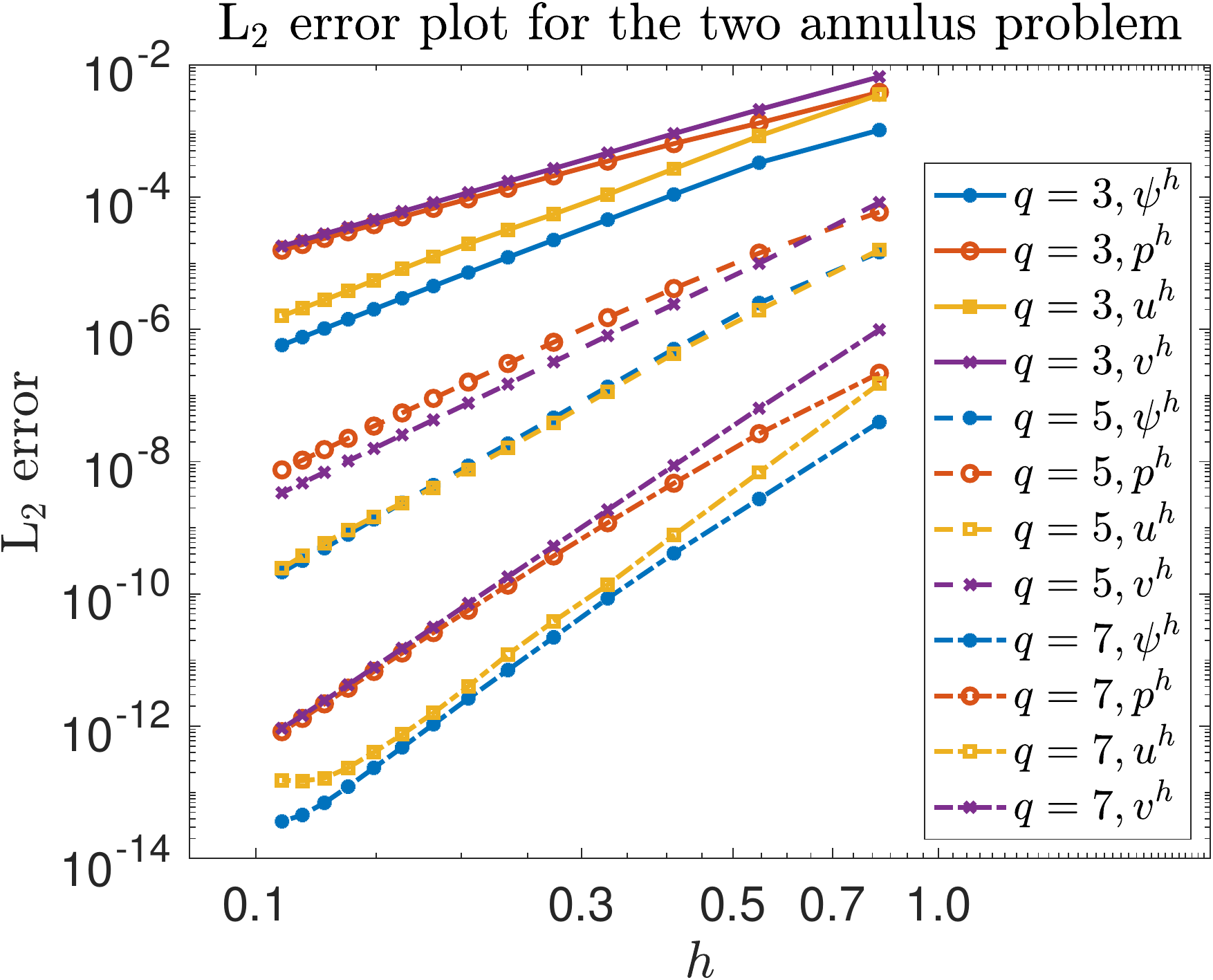}
\caption{Error plot for the two annulus example. \label{err_plot_annulus}}
\end{center}
\end{figure}

We evolve the solution until the time $t = 1$ and plot the L$_2$ error at the final time in Figure \ref{err_plot_annulus}. The convergence rates are shown in Table \ref{table_annulus}. For $q = 3$ and 5 we observe optimal rates, i.e $q+1$ for $\psi^h$ and ${\bf u}^h$ and $q$ for $p^h$ and ${\bf v}^h$. For $q = 7$ we still observe optimal rate of convergence for ${\bf v}^h$ and close to optimal rates for ${p}^h$ and $\psi^h$ while  the rate of convergence for ${\bf u}^h$ only is $q$. 
\begin{table}
\centering
\caption{Computed rates of convergence with $q_\psi=q_u=q$ and $q_p=q_v=q-1$ for the two-annulus example.\label{table_annulus} }
\begin{tabular}{c c c c c}
\toprule
$q$ & $\psi^h$ & $p^h$ & $\bs{u}^h$  & $\bs{v}^h$  \\
\midrule
3 & 4.03 & 2.87 & 3.98 & 2.96 \\
5 & 5.92 & 4.92 & 5.78 & 4.97 \\
7 & 7.57 & 6.70 & 6.98 & 6.98 \\
\bottomrule
\end{tabular}
\end{table} 
The rates are computed using the 5 finest grids. The exception is for $u^{h}$ and $q = 7$ where we exclude the four finest grids due to finite precision effects and use the next five to compute the rate of convergence.   

\subsection{Source location via time reversal}
In this example we consider an example of source location via time-reversal in a coupled fluid solid problem. Time Reversal (TR) is a technique to solve certain inverse problems such such as source localization. The key ingredient in TR is the fact that the wave equation (containing only second order time derivatives) is invariant to the change of variables $t = -t$. This means that when evolving a Cauchy problem governed by the  wave equation until time $T$, starting from some localized initial data, it is possible to recover the initial data exactly by simply evolving the solution at time $T$ backwards in time to time $t=0$. More surprisingly, the localization of the initial data persists to a high degree even when the data is scarce, for example when recorded for a finite time and at a finite number of recorders. Time reversal has received considerable attention in recent years, see for example the review by Givoli \cite{TR_Givoli}. 

Here we consider an example of a fluid on top of a solid inside the domain $(x_1,x_2) \in  [-1,1] \times [-2,2]$ and separated by an interface described by the curve 
\[
\gamma(x_1,x_2) = (0,0.025 \sin ( n \pi x_1)).
\]
The boundary condition on the top and bottom are homogenous Dirichlet conditions and on the sides we impose homogenous traction and Neumann conditions in the solid and fluid respectively. We solve until time 4 and use $70 \times 70$ elements with degree 5 in the fluid and the solid. The initial data representing the source is confined to the solid and is centered at $(0,-0.4)$. We consider two cases, either the pressure wave   
\[
u_1 = x_1 f(r), \ \ u_2 = (x_2+0.4) f(r),
\]
or the shear wave
\[
u_1 = -(x_2+0.4) f(r), \ \ u_2 = x_1 f(r).
\]
Here 
\[
f(r) = e^{-\frac{36}{0.4^2} (x_1^2+(x_2+0.4)^2)}.
\]
The fluid is assumed to be at rest and the velocity in the solid is also taken to be zero. 
\begin{figure}[htbp]
\begin{center}
 \includegraphics[width=0.19\textwidth]{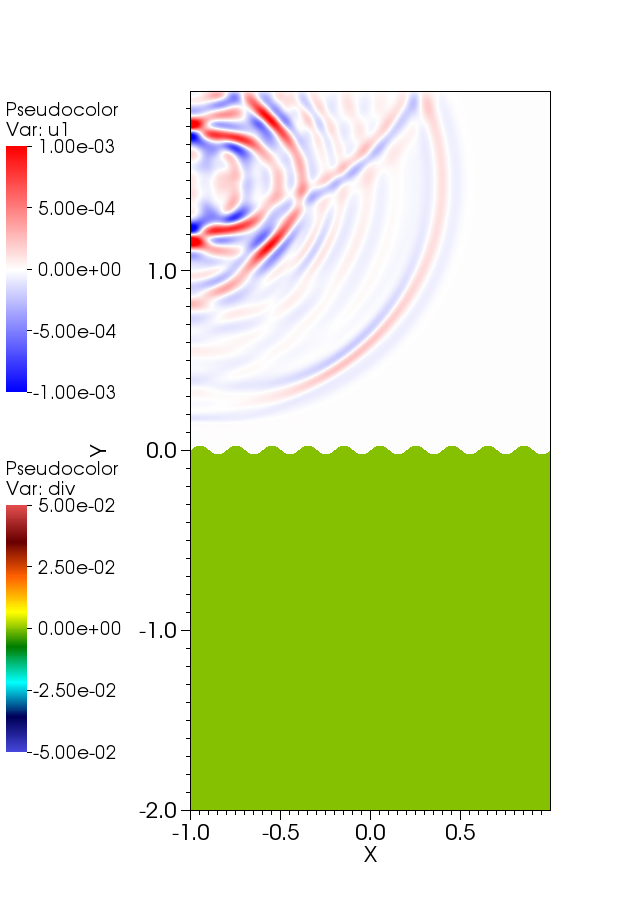}
 \includegraphics[width=0.19\textwidth]{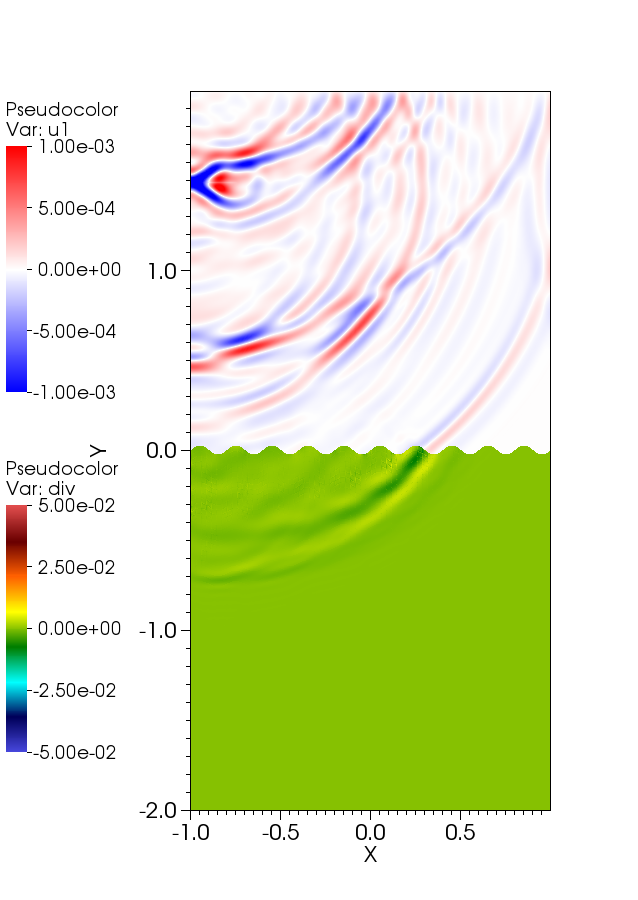}
 \includegraphics[width=0.19\textwidth]{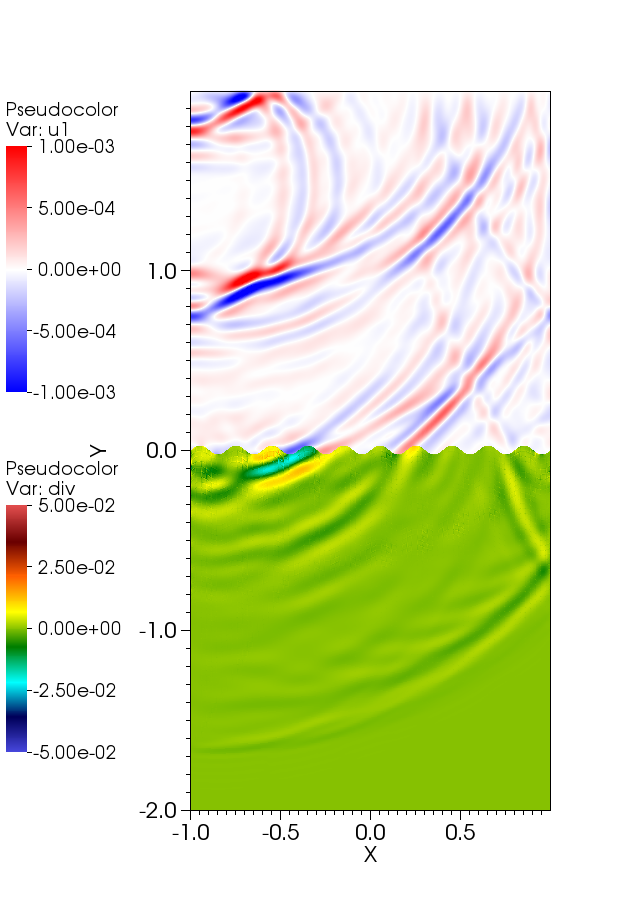}
 \includegraphics[width=0.19\textwidth]{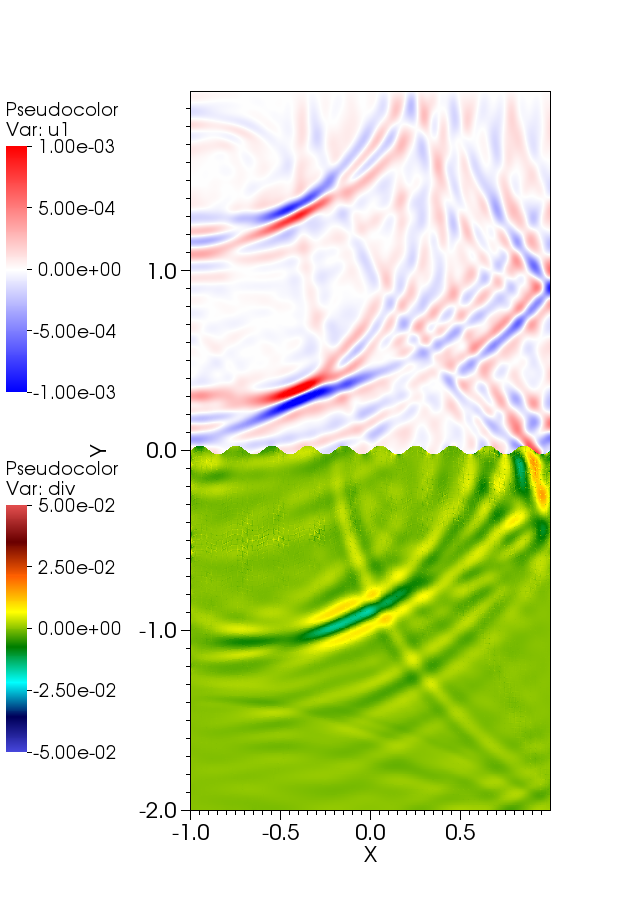}
 \includegraphics[width=0.19\textwidth]{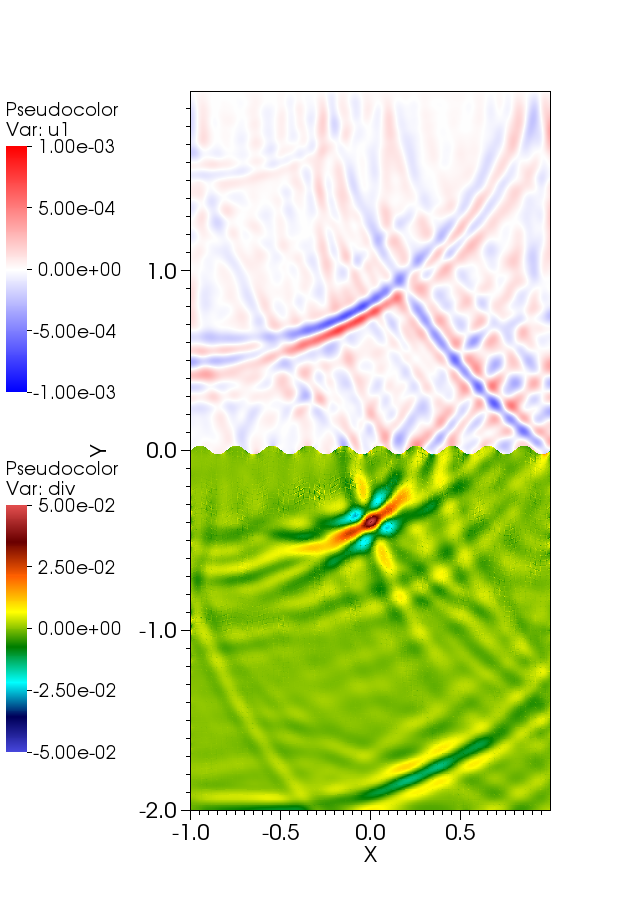}
 \caption{Backward propagation of data recorded due to the pressure wave.   \label{fig:tr3}}
\end{center}
\end{figure}

This problem illustrates that TR has ability to distinguish between an explosive source of compressional and source generating a shear wave from underwater measurements, something that is of interest, e.g. for monitoring of tests of nuclear bombs. 

During the forward simulation we record the velocity potential in the fluid along the straight line connecting the points (0.0286,1.5714) and (0.2571,1.5714). The recorded signal is then re-transmitted. As can be seen in Figure \ref{fig:tr3} the reverse signal is re-focused at the location of the source at the initial time.    

\begin{figure}[htbp]
\begin{center}
 \includegraphics[width=0.3\textwidth]{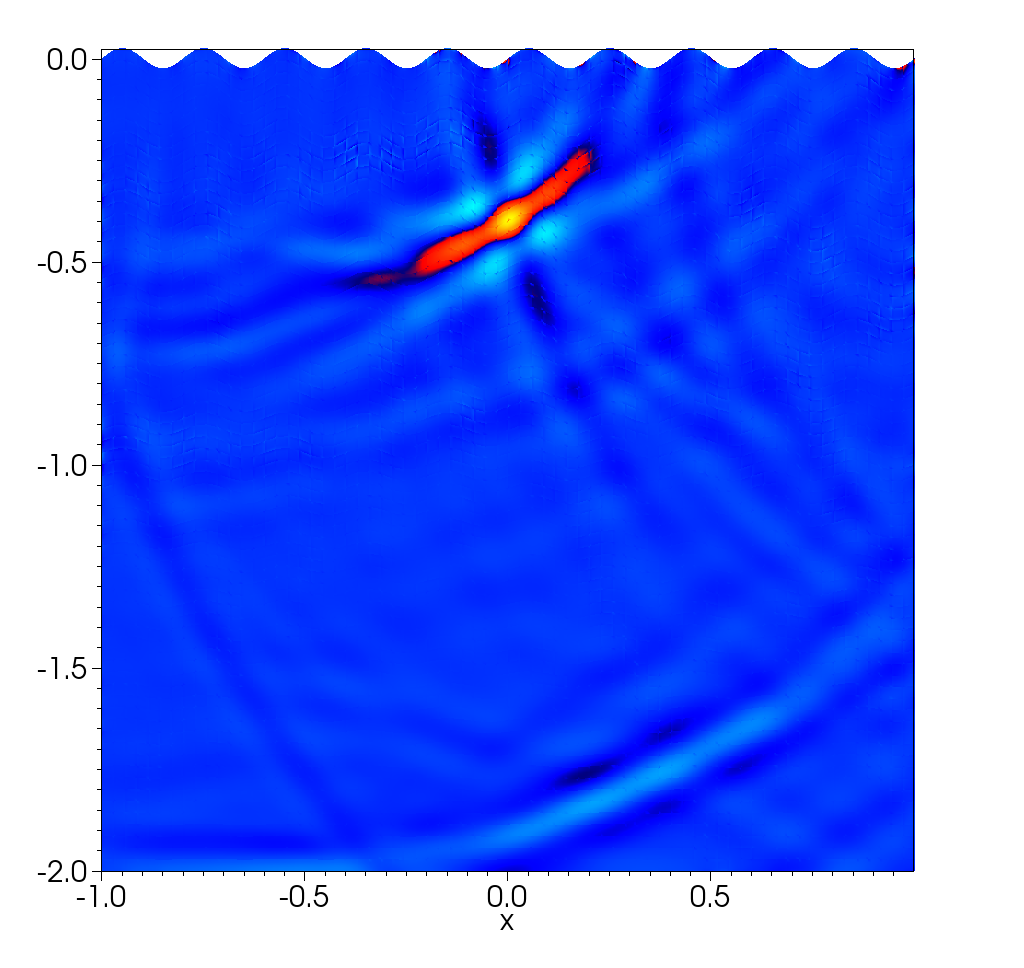}
 \includegraphics[width=0.3\textwidth]{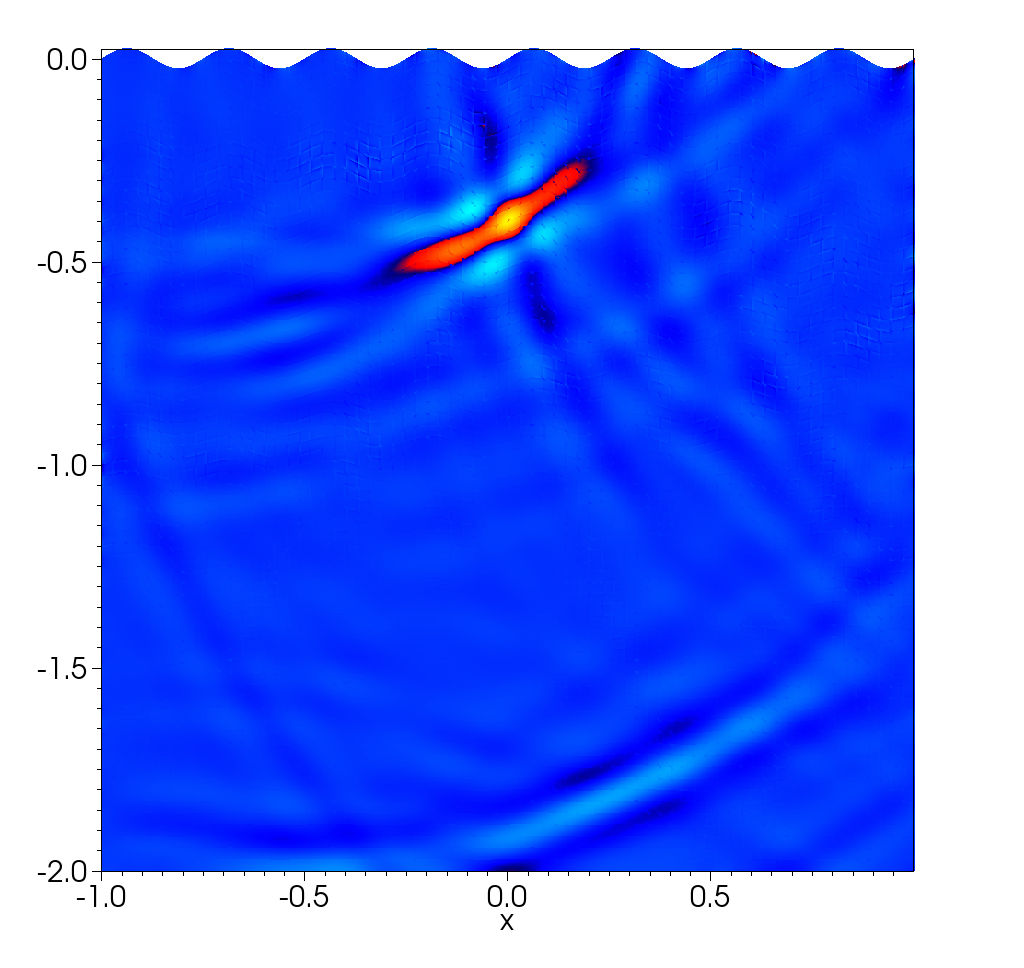}
 \includegraphics[width=0.3\textwidth]{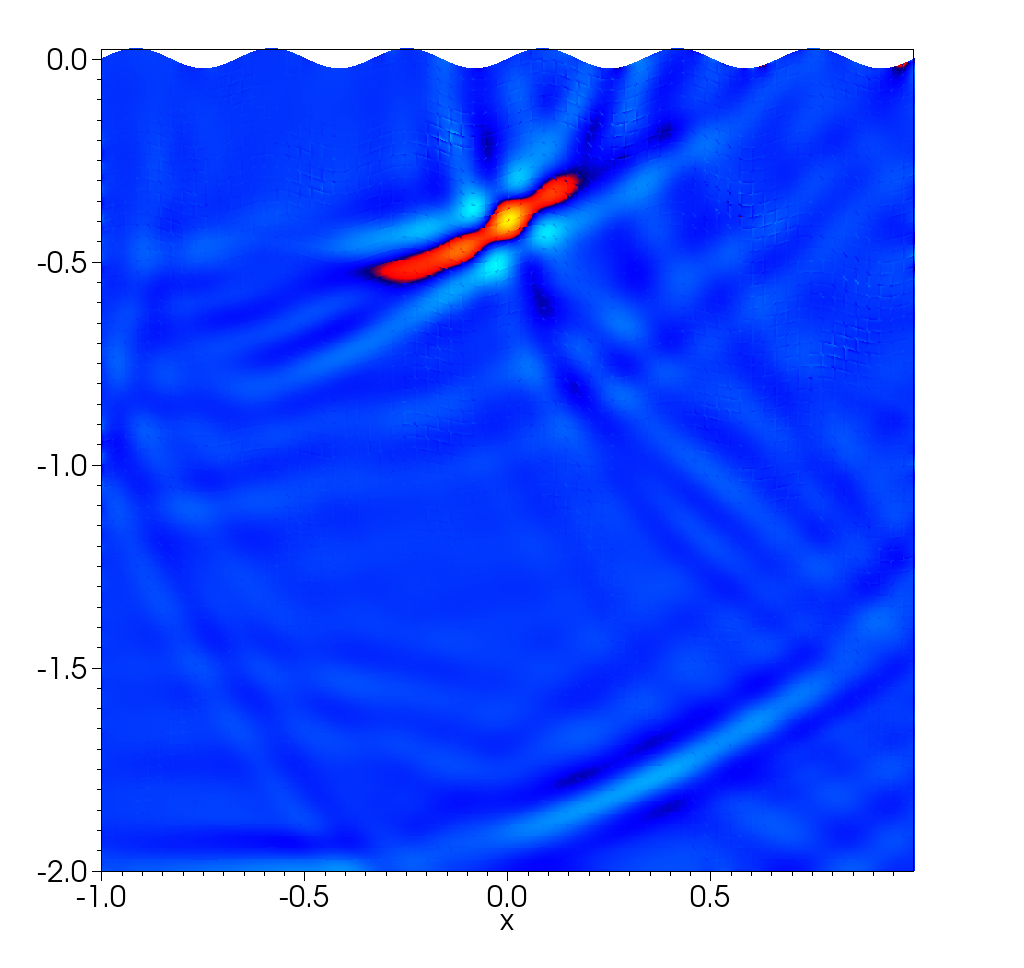}
 \includegraphics[width=0.3\textwidth]{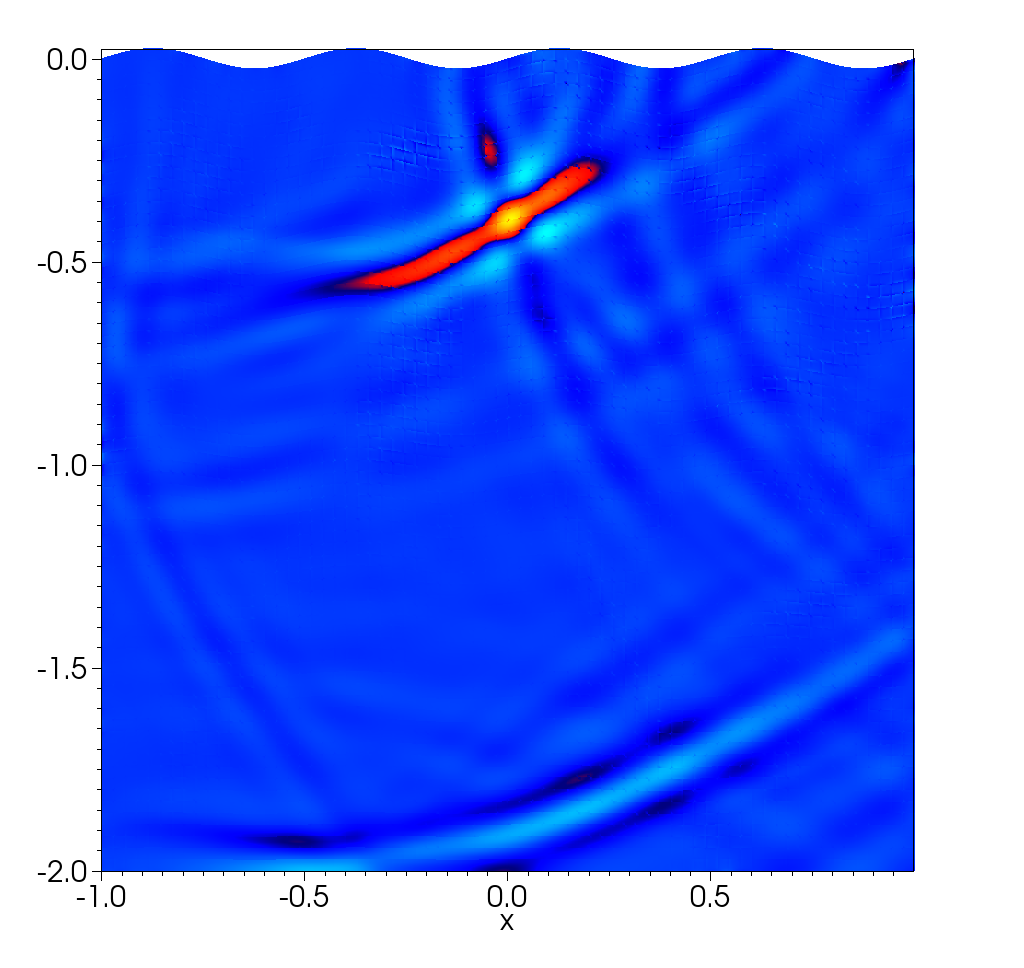}
 \includegraphics[width=0.3\textwidth]{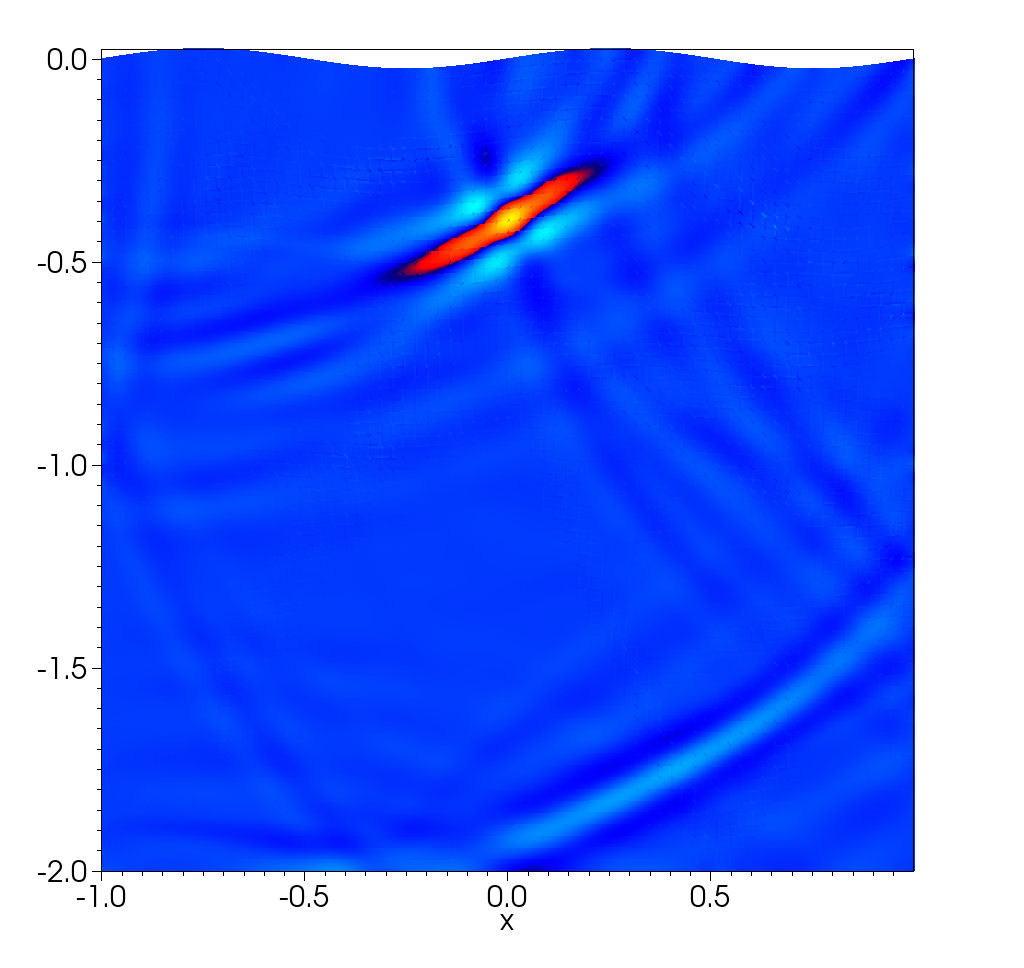}
 \includegraphics[width=0.3\textwidth]{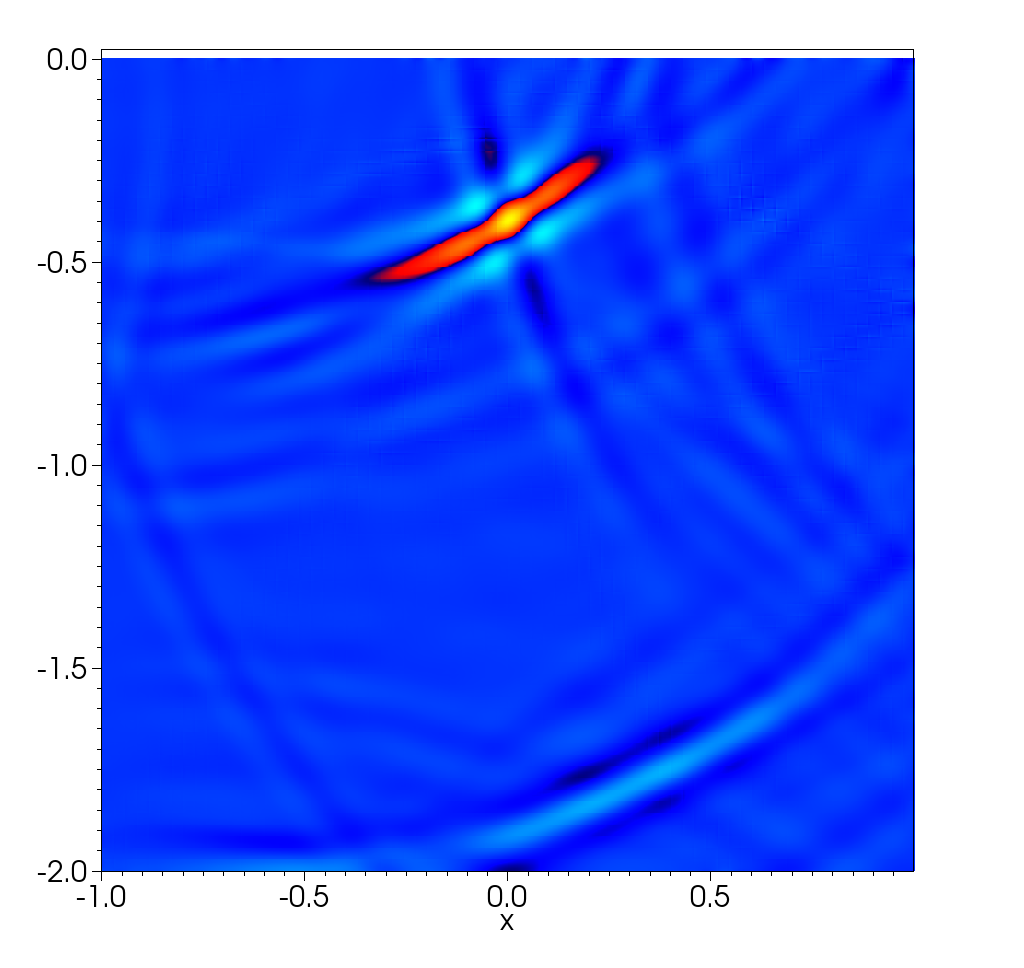}
 \caption{The divergence of the displacement at the final time. Pressure wave initial data.  \label{fig:tr1}}
\end{center}
\end{figure}

\begin{figure}[htbp]
\begin{center}
 \includegraphics[width=0.3\textwidth]{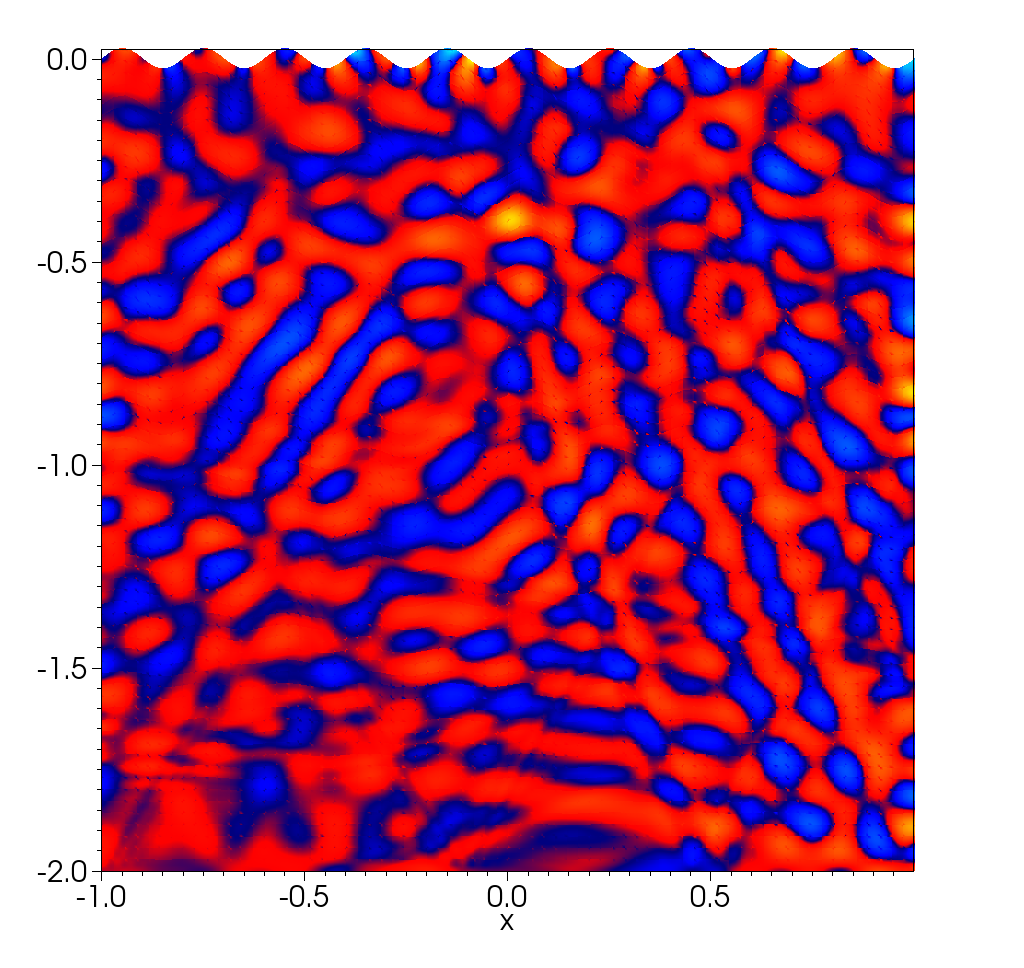}
 \includegraphics[width=0.3\textwidth]{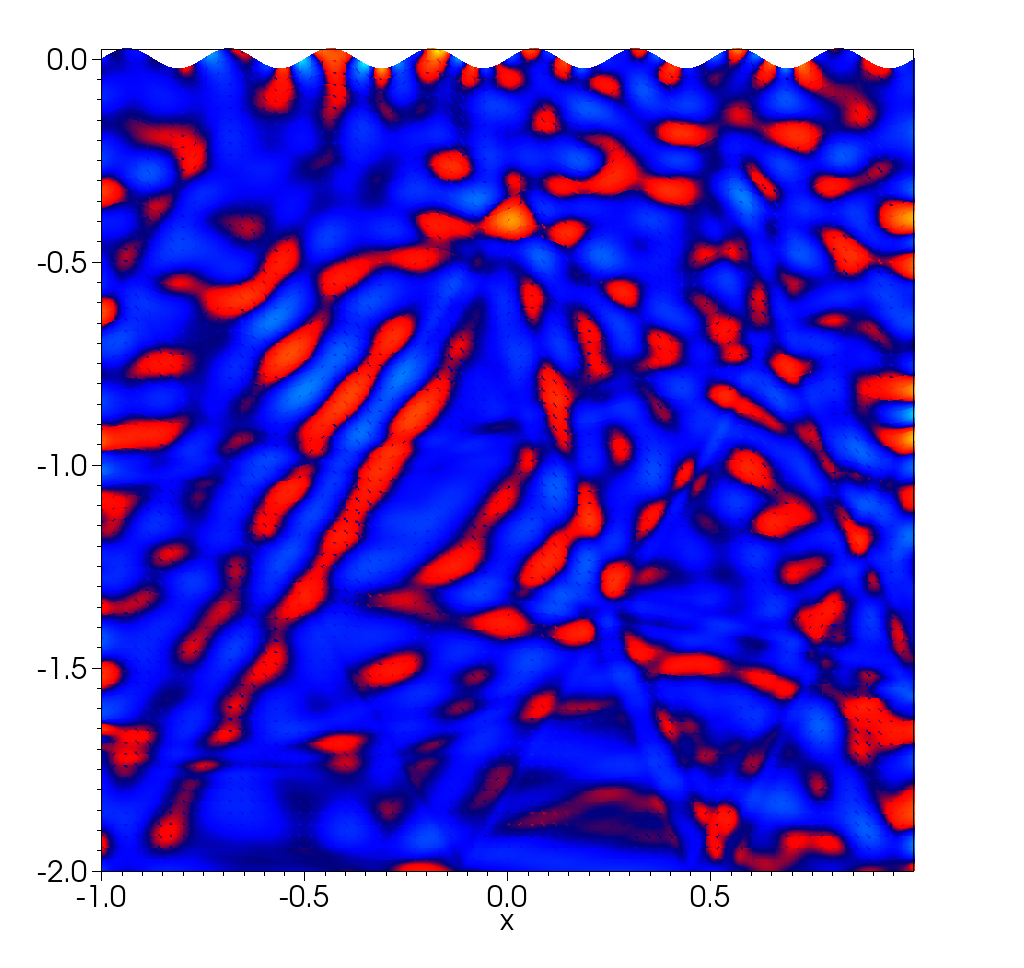}
 \includegraphics[width=0.3\textwidth]{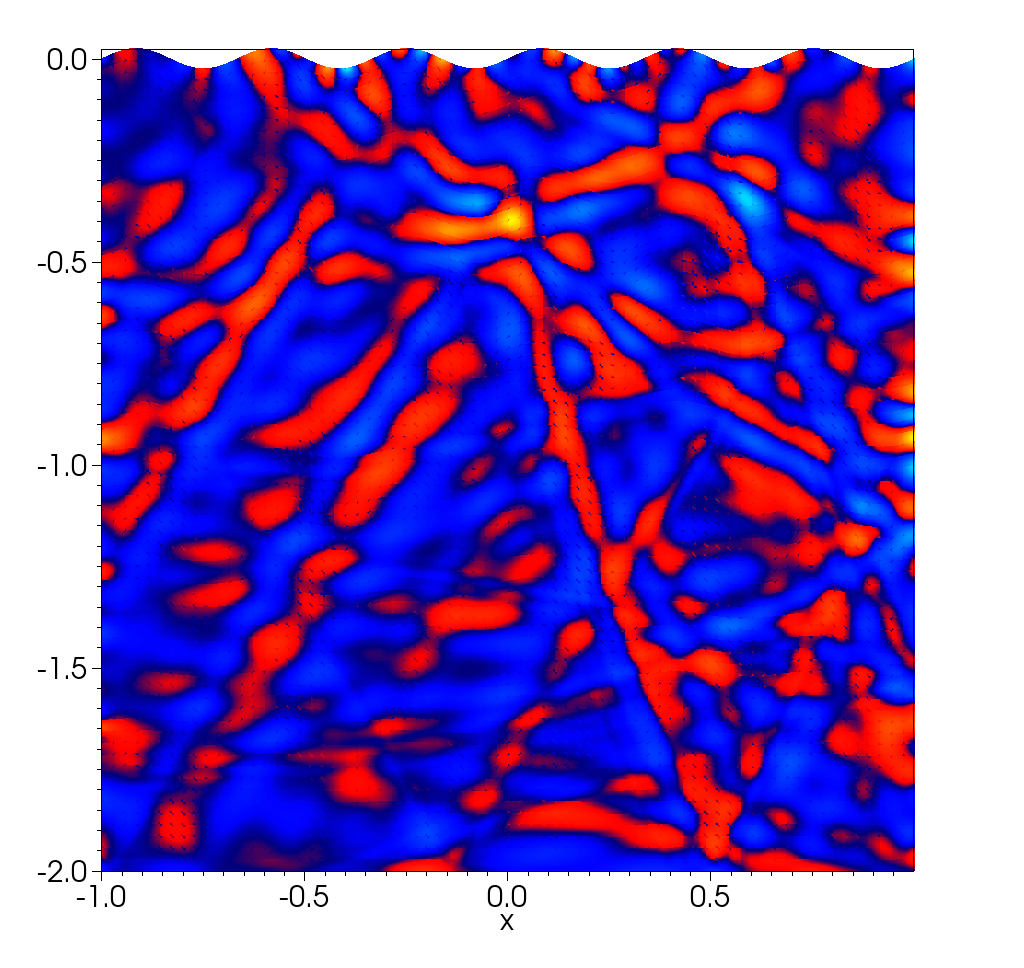}
 \includegraphics[width=0.3\textwidth]{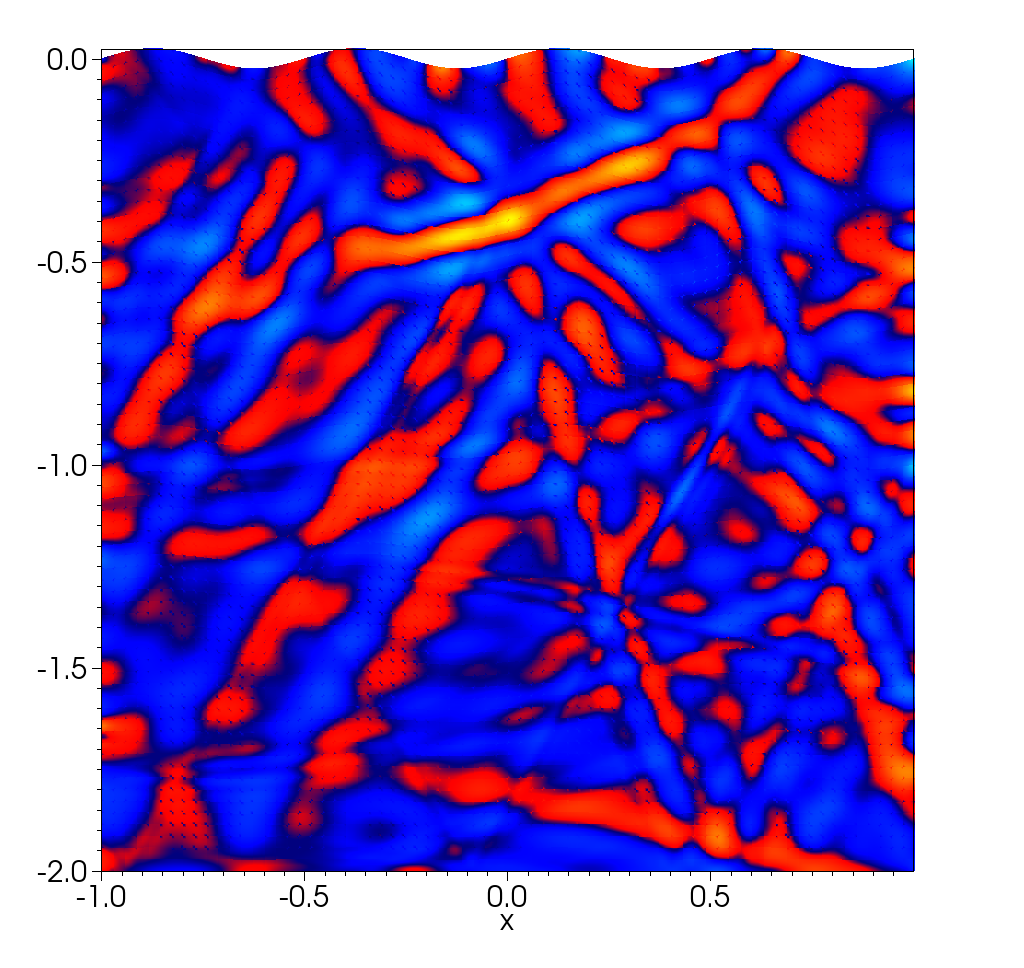}
 \includegraphics[width=0.3\textwidth]{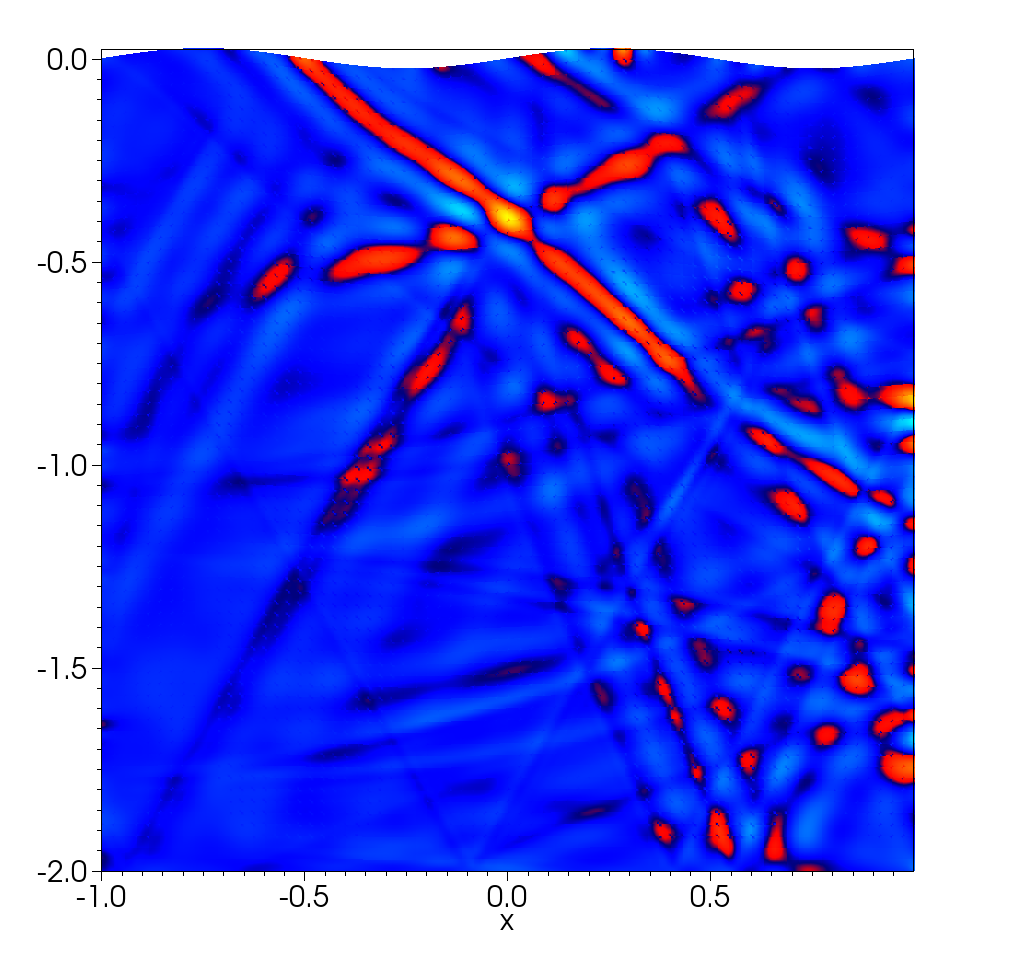}
 \includegraphics[width=0.3\textwidth]{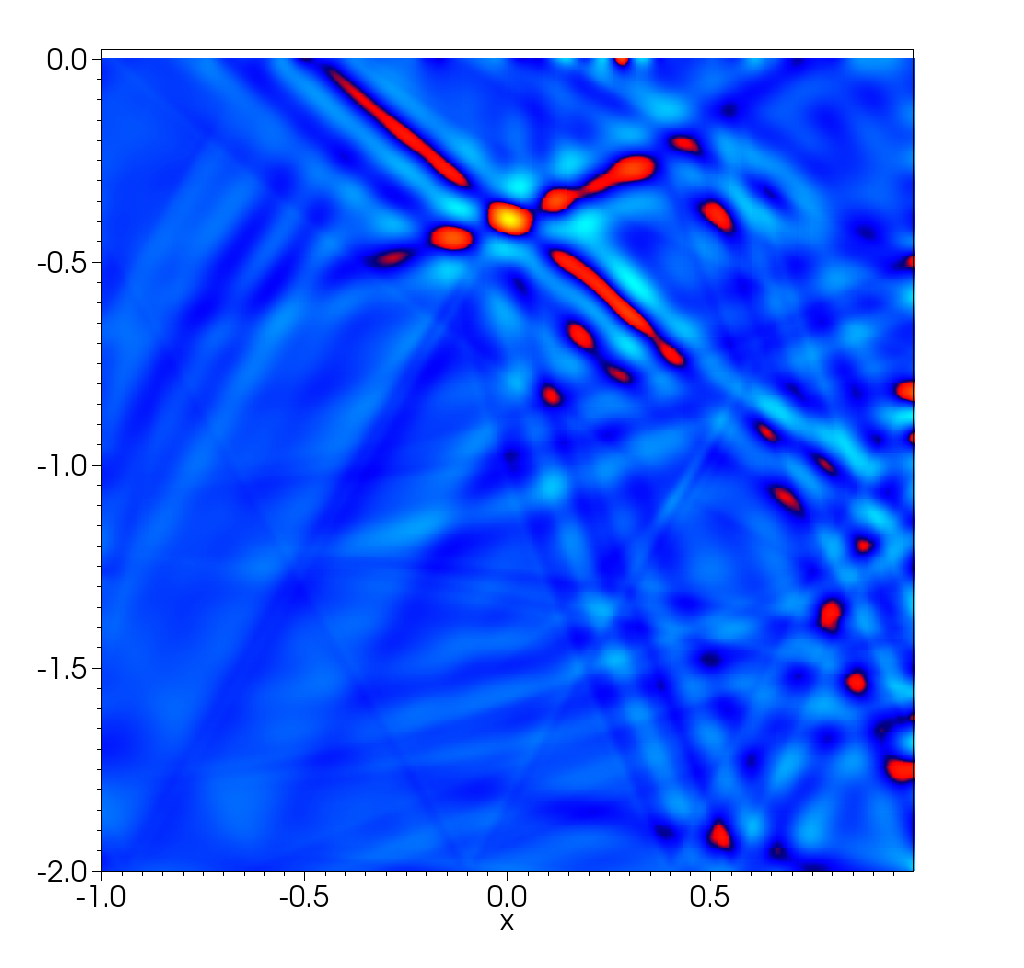}
 \caption{The rotation of the displacement at the final time. Shear wave initial data. \label{fig:tr2}}
\end{center}
\end{figure}

In Figure \ref{fig:tr1} and \ref{fig:tr2} we display the divergence and the rotation of the displacement at the final time for interfaces with $n = 10,8,\ldots,0$. The initial data for Figure \ref{fig:tr1} was the pressure wave and the shear wave for Figure \ref{fig:tr2}. The color scheme was picked to emphasize contrasts but is the same for both sets of initial data. 

For the case of the pressure wave we see a rather distinct focus for all of the interfaces with perhaps a slight improvement with increasing structure, see Figure \ref{fig:tr1}. For the case of a shear wave we do not get a clear focus for any of the interfaces and here the focusing appears to deteriorate with increased structure of the interface.   


\subsection{Inverse Problems}
In this section we present some applications of our method applied to inverse problems. 
\subsubsection{Interface inversion}
\begin{figure}[h]
\begin{center}
 \includegraphics[width=0.35\textwidth]{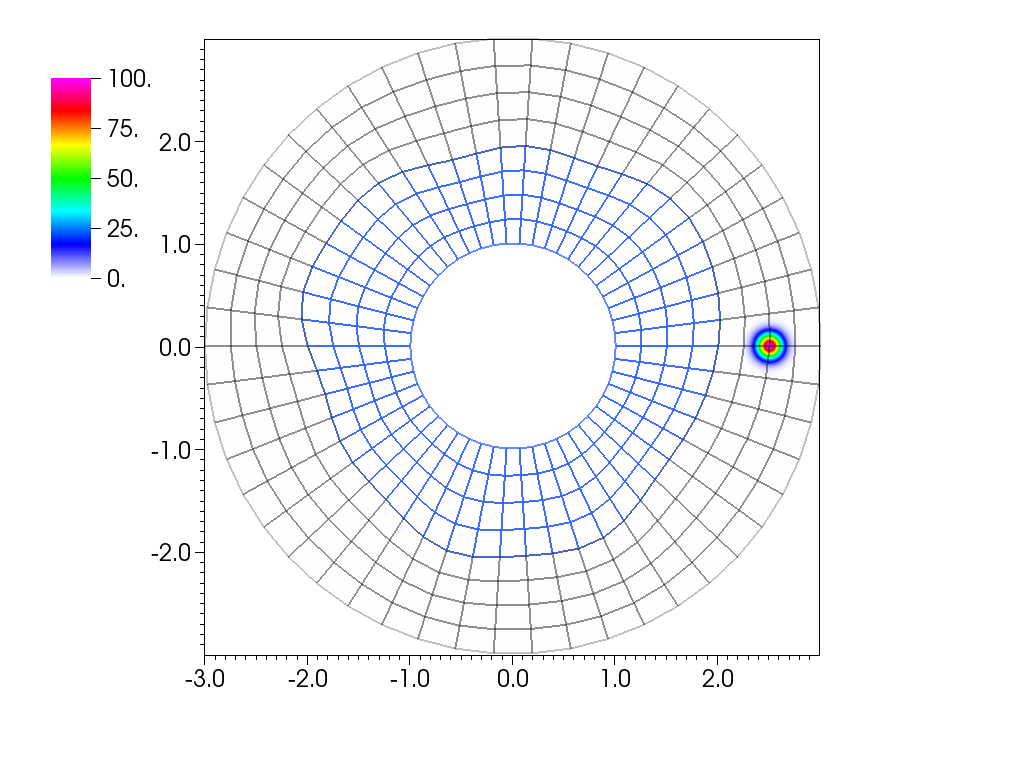}
 \caption{The geometry and initial data used in the interface inversion problem. \label{fig:surf1}}
\end{center}
\end{figure}
In this example we consider the inverse problem of finding the shape of an internal interface between an outer annular shaped region occupied by a fluid and an inner annular shaped region occupied by a solid. To do this we create synthetic data by recording the velocity potential $\psi$ at 50 equidistant locations on the surface of the fluid (where we impose a homogenous Neumann condition). The inner and outer radii are 1 and 3 respectively and we are considering interfaces 
described by the periodic function 
\begin{equation*}
(x_{1, \rm I}(\theta), x_{2, \rm I}(\theta)) =  ( (1+\delta r (\theta)) \cos(\theta), (1+\delta r (\theta)) \sin(\theta)), \ \ \delta r (\theta) = \sum_{k = 1}^{8} A_k \sin ((k+1)\theta), \ \ \theta \in [0, 2 \pi],
\end{equation*}
and the synthetic data is obtained using the coefficients $A_k, k=1,\ldots,8$ 
\begin{equation}
[0.002, 0.050, -0.001, 0.008, -0.003, -0.006, -0.010, 0.010].
\end{equation}
The initial fields are all taken to at rest except for  
\[
\psi(x_1,x_2,0) =100 e^{-72((x_1-2.5)^2+x_2^2))}, \ \ (x_1,x_2) \in \Omega_{f}.
\]
The grid used and the initial data are displayed in Figure \ref{fig:surf1}.

The speed of sound in the homogenous fluid is taken to be one and the lame parameters in the solid are $\lambda = 2$ and $\mu = 0.2$ and the density is $\rho=1$. The solution is recorded until the final time is 25.  

\begin{figure}[h]
\begin{center}
 \includegraphics[width=0.32\textwidth]{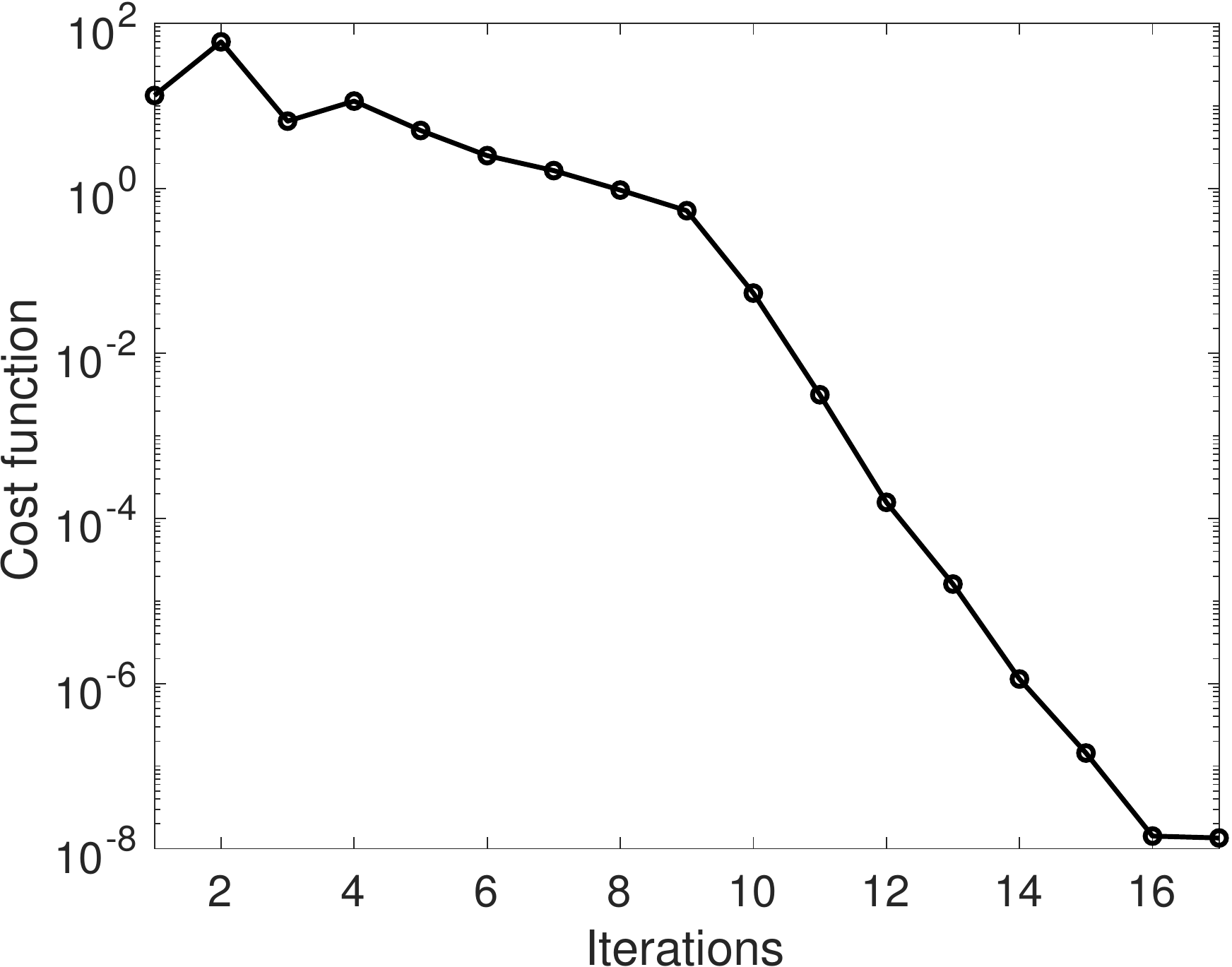} \ \ 
 \includegraphics[width=0.31\textwidth]{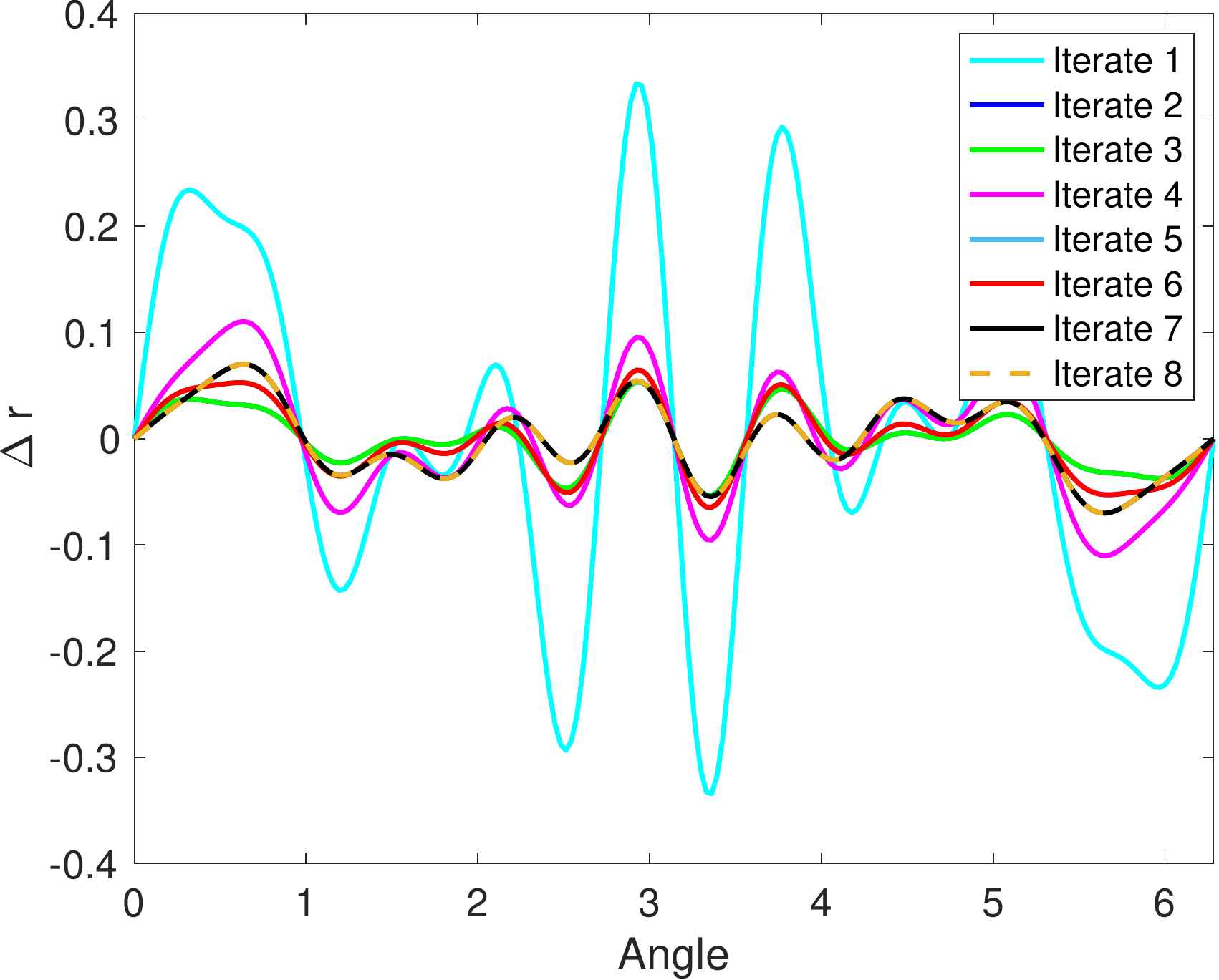} \ \ 
 \includegraphics[width=0.31\textwidth]{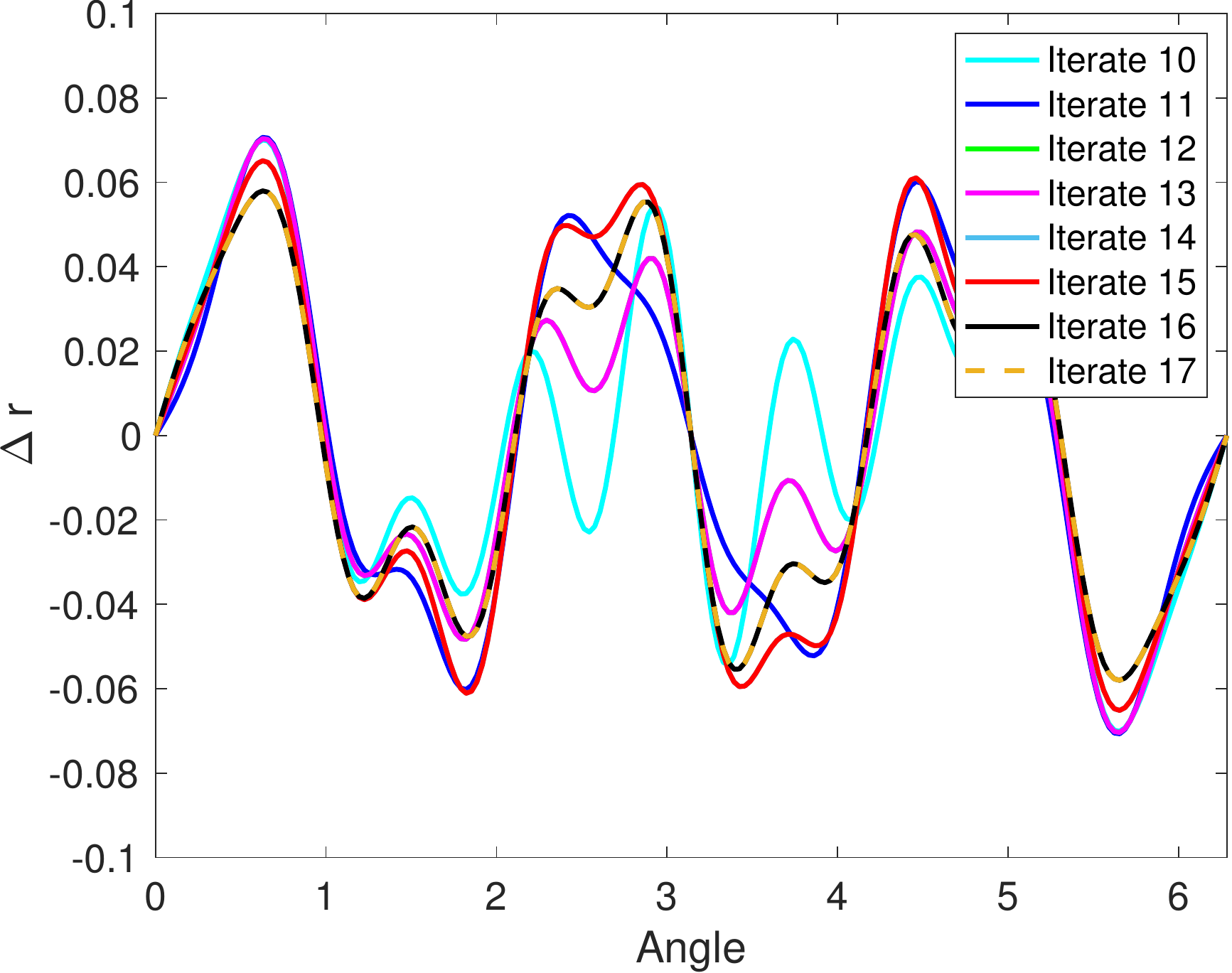}
 \caption{To the left the the reduction of the cost functional as a function of the number of iterations. The middle and right figures display the shape of the interface for different iterations.  \label{fig:surf2}}
\end{center}
\end{figure}

The cost (misfit) functional that we minimize is the sum of the squares of the $L_2$ norms in time of the difference between the synthetic data and the simulated data. During the minimization we impose the bounds $| A_k | \le 0.1$. To recover $A_k, k=1,\ldots,8$ 
we use the \verb+L-BFGS-B+ algorithm by Byrd et al. \cite{L-BFGS} and since the number of variables is small we simply use forward differences to to approximate gradients.  
The reduction of the cost functional and the change in the interface are displayed in Figure \ref{fig:surf2}. As can be seen the convergence is quite rapid and a reduction of the cost function by ten orders of magnitude only takes about 15 iterations. We note that for the first few iterates it is necessary to enforce the bounds on $A_k$.

\subsubsection{A simple material model inversion}
In this example we perform a full wave inversion of the compression wave speeds in the solid. The fluid domain is $\Omega_{f} \in [-1,1] \times [0,2]$ and the domain of the solid is $\Omega_{s} \in [-1,1] \times [-2,0]$. Each of the domains is discretized by $5 \times 5$ square elements. We impose homogenous Dirichlet boundary conditions on the vertical sides, a homogenous Neumann boundary condition on top of the fluid domain and we set the bottom surface in the solid to be free of traction. 

In the fluid the speed of sound is one, and in the solid we have in each of the elements $\mu = 2$ and $\lambda = 4 + \delta \lambda_i, \, i = 1,\ldots, 25$, where $\delta \lambda_i \in [0, 1]$. 

\begin{figure}[htbp]
\begin{center}
\includegraphics[width=0.48\textwidth]{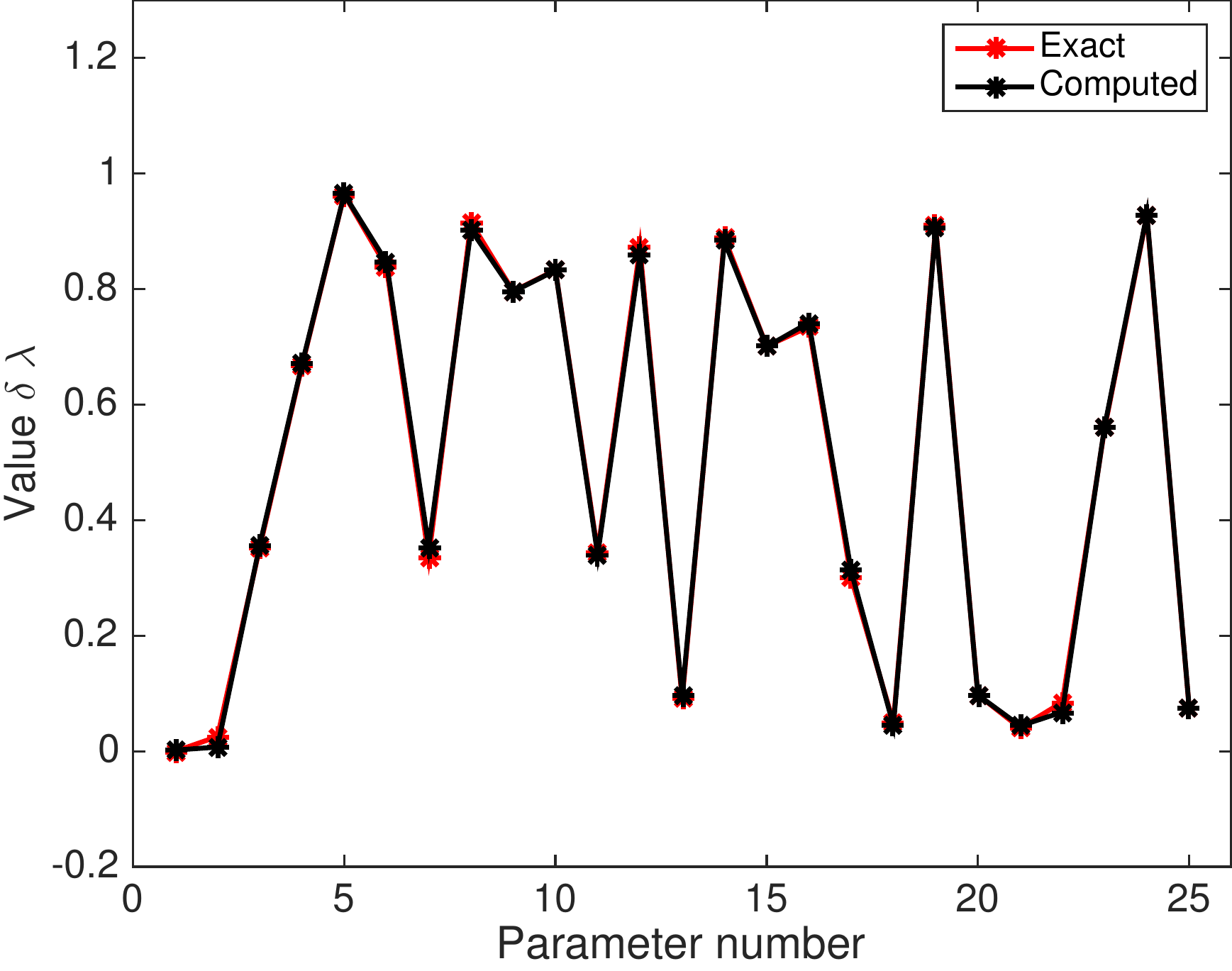} \ \
\includegraphics[width=0.465\textwidth]{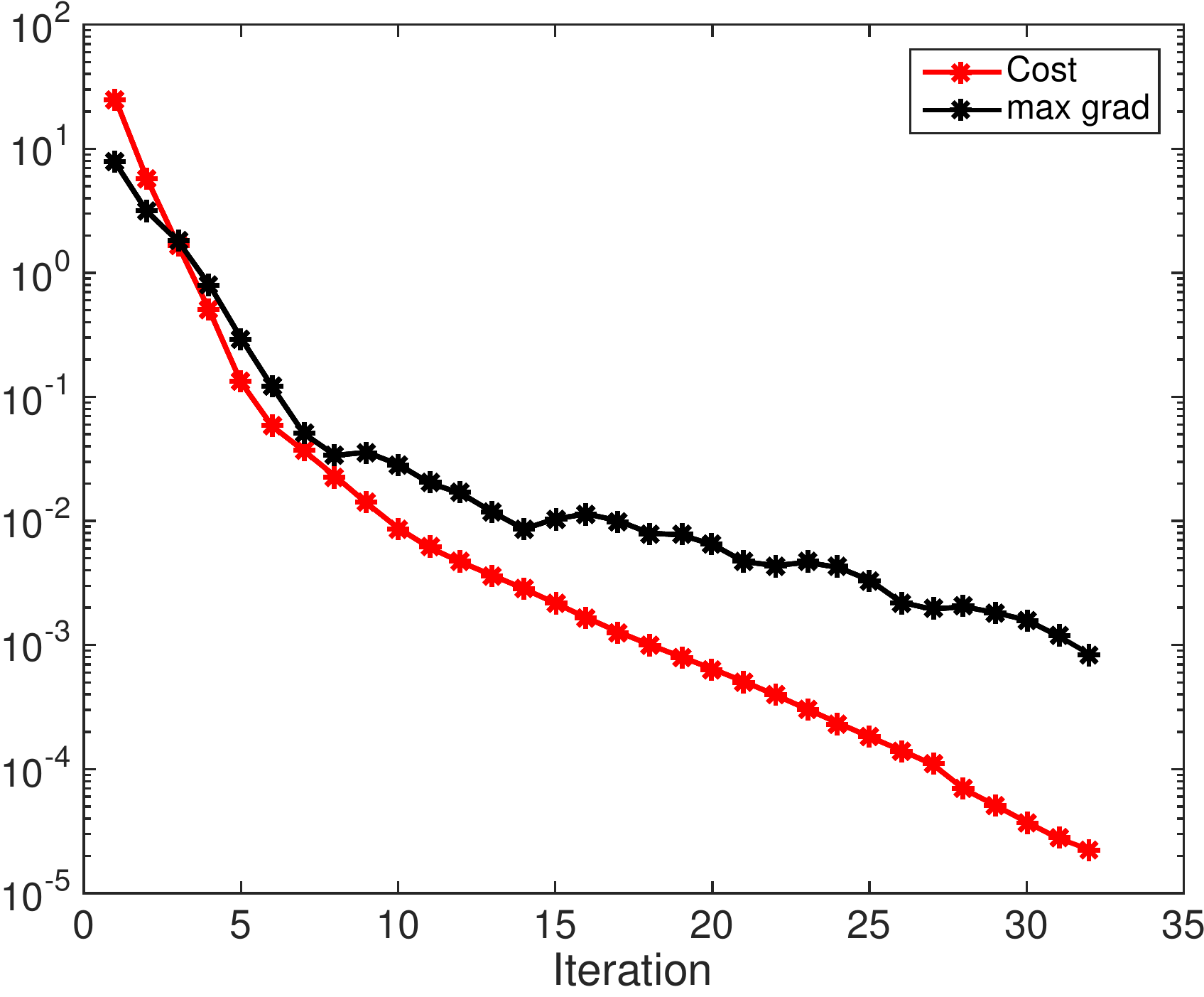}
\caption{To the left we display the values of the perturbations of $\lambda$ for the material (in red) and the inverted material (in black). To the right we display the cost function and the maximum norm of the gradient as a function of the number of iterations. \label{fwe:1}}
\end{center}
\end{figure}

The solid and fluid are initially at rest and the solid is forced by a point source in the point $(x_1,x_2) = (0.1,1.8)$. The amplitude of the source is 
\[
A(t) =  -2w_0^2 (t-t_0) e^{-(w_0(t-t_0))^2}.
\]
We set in the experiment $t_0 = 1$ and $w_0 = 6$, and record the solution by 9 receivers placed at the surface of the fluid starting at $x_1 = -0.8$ and with a spacing of 0.2 from time 0 to 10. 

The forward solve uses the upwind flux everywhere and polynomials of degree 8. The solution is advanced until time 8 by a Taylor series method with 8 time derivatives. The time step is set by the stability restrictions in the solid. Precisely we have $dt/h =0.15/(\sqrt{2\mu+\lambda} (q+1.5))$ with $\mu = 2$ and $\lambda = 4$.   

\begin{figure}[htbp]
\begin{center}
 \includegraphics[width=0.3\textwidth]{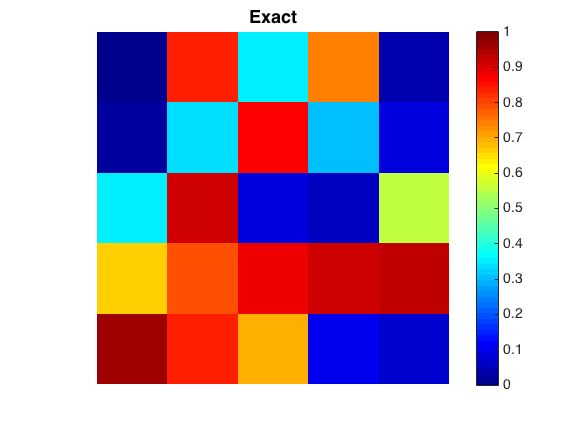}
 \includegraphics[width=0.3\textwidth]{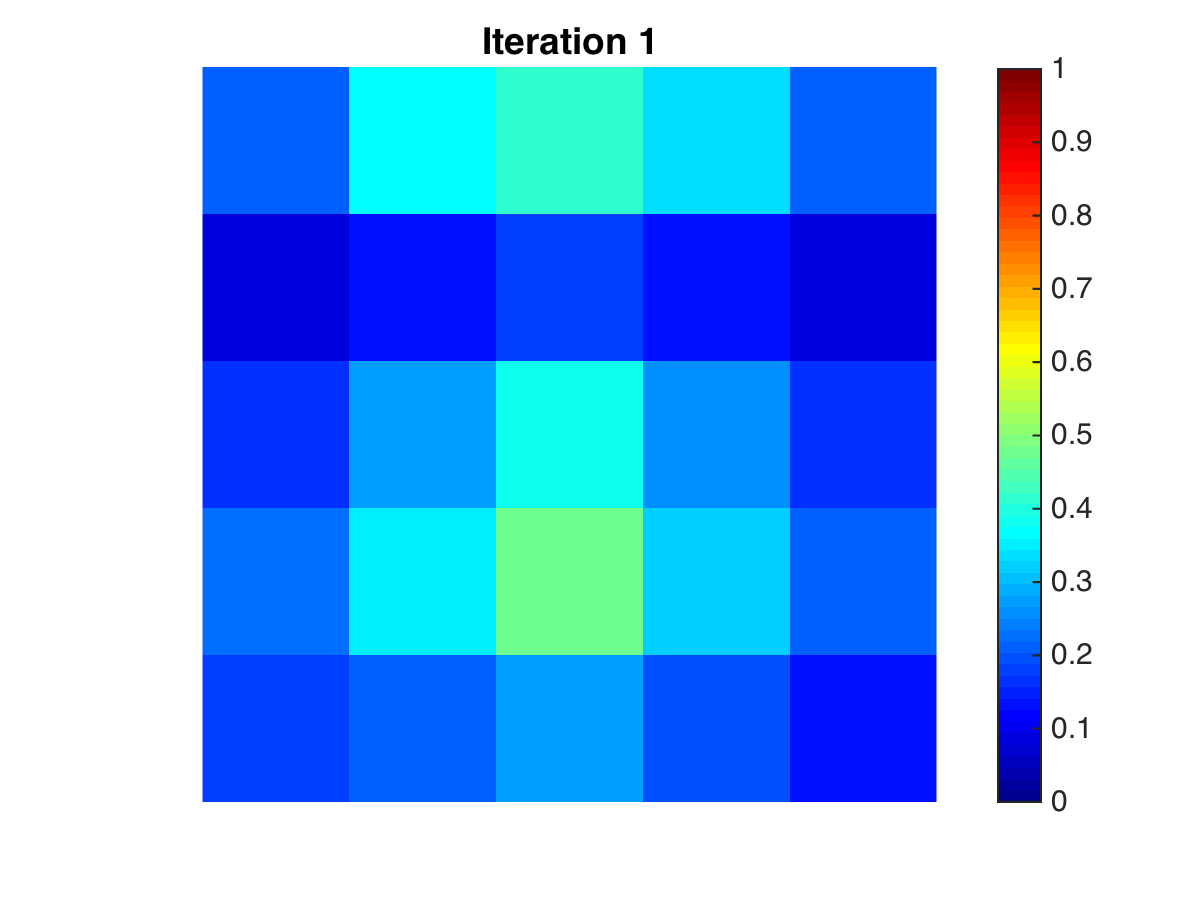}
 \includegraphics[width=0.3\textwidth]{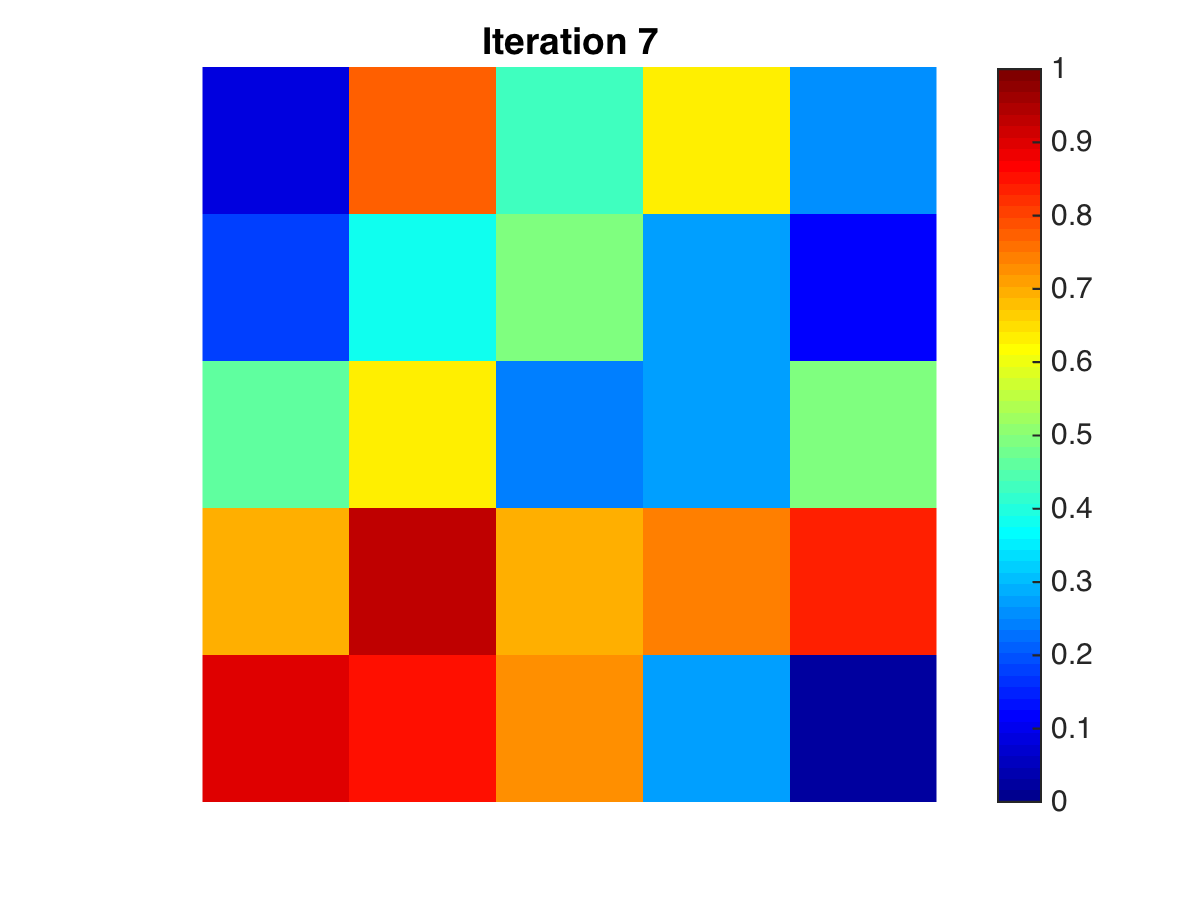}
 \includegraphics[width=0.3\textwidth]{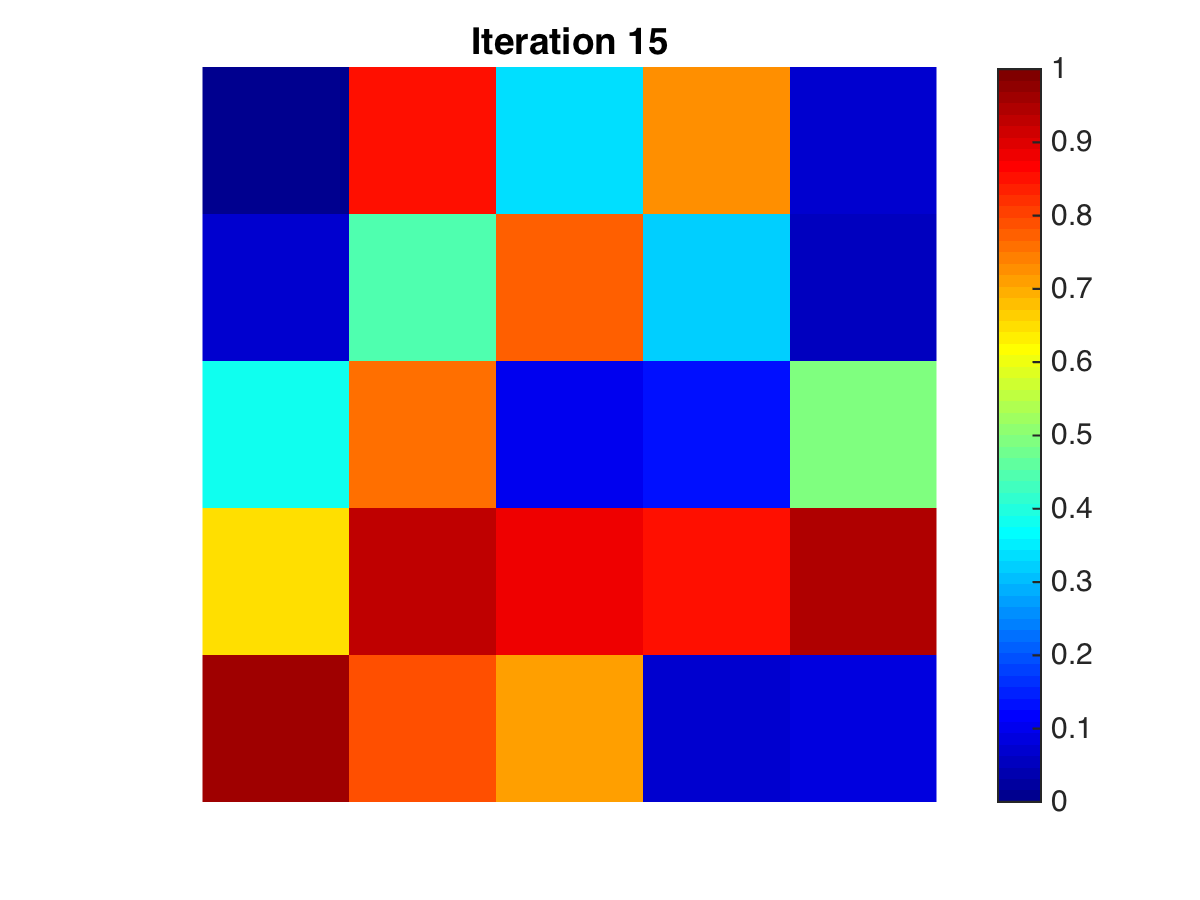}
 \includegraphics[width=0.3\textwidth]{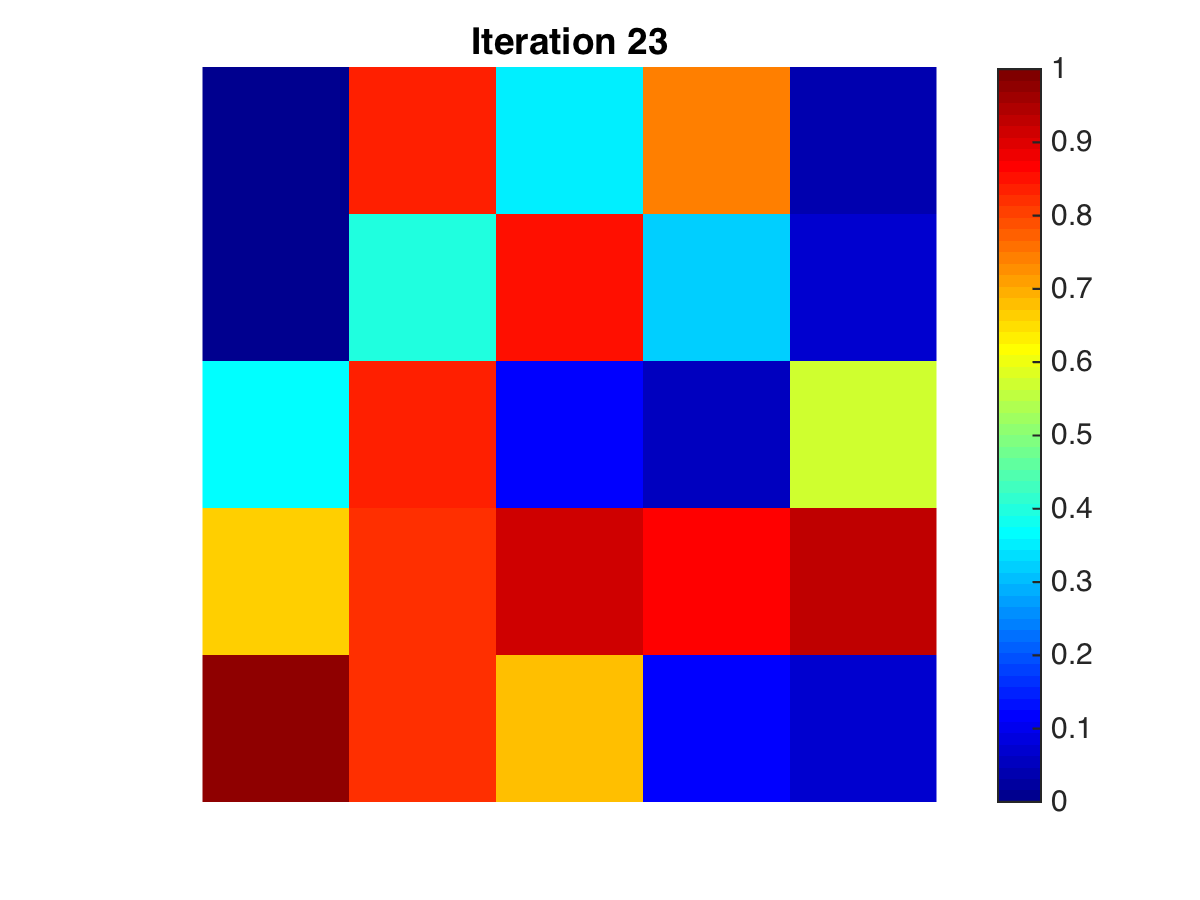}
 \includegraphics[width=0.3\textwidth]{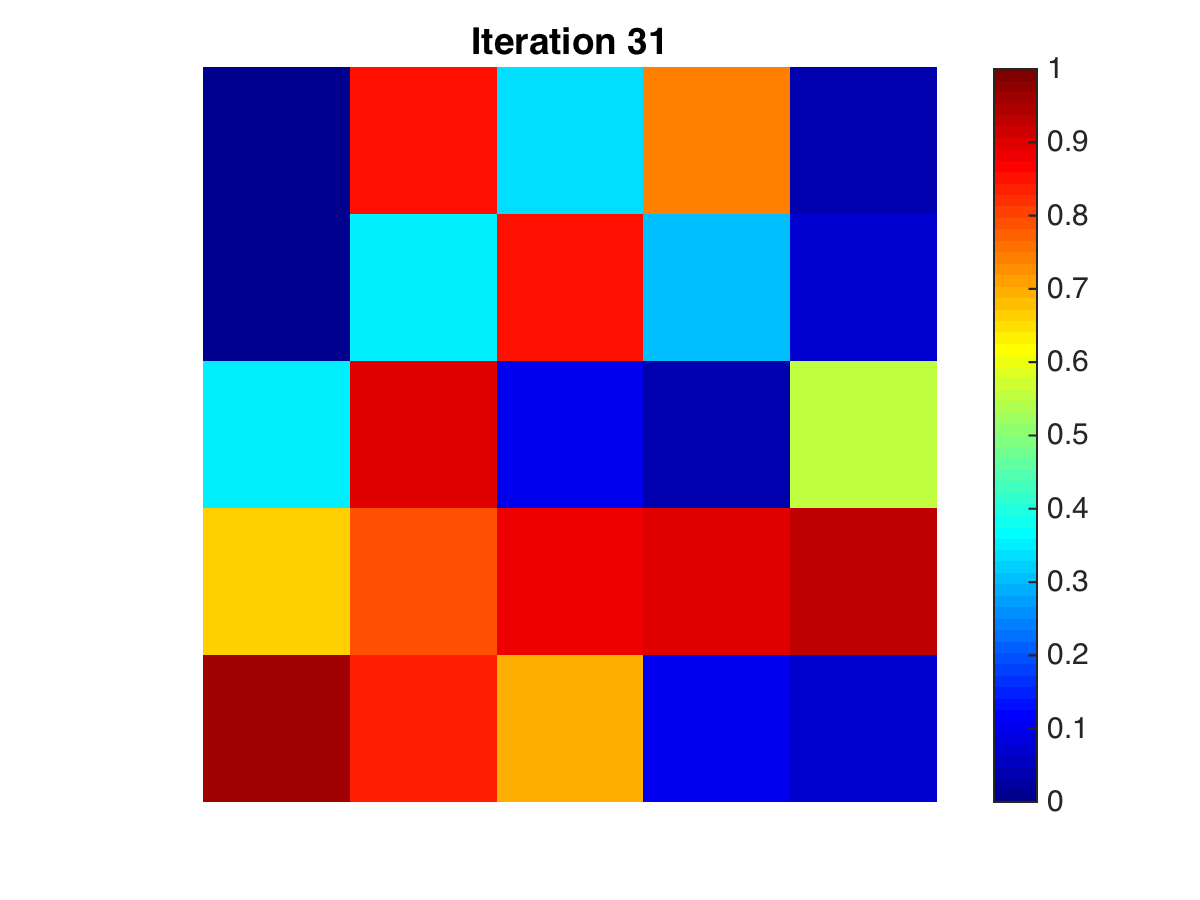}
 \caption{Progression of the inversion process. Plotted is the material model with increasing number of iterations in the minimization process. \label{fwe:2}}
\end{center}
\end{figure}

We invert for the perturbation to the material parameters, $\delta \lambda_i$ by minimizing the sum of the misfits (measured in the $L_2$-norm in time) in the nine receivers. The exact data is obtained by a single forward solve with the known material parameters. To minimize the total misfit we use the Broyden-Fletcher-Goldfarb-Shanno algorithm, \cite{Nocedal2006NO}, with a backtracking line search utilizing the Armijo-Goldstein condition. The initial guess for the perturbations is all zeros. Here we form the gradients by a finite difference approximation. The results, displayed in Figures \ref{fwe:1} and \ref{fwe:2} illustrate the convergence of the minimization process.  

\section{Conclusion and future work}\label{sec-conclusion}
We have studied a wave propagation problem in a coupled fluid solid region, with a focus on the fluid solid interface. Wave propagation is modelled by the wave equation in terms of velocity potential in the fluid, and the elastic wave equation in displacement in the solid. The energy based discontinuous Galerkin method is used to discretize the governing equations in space. We have derived both energy conserving and upwind numerical fluxes. Our test problems show that upwind fluxes perform best, with either central or alternating coupling. Numerical experiments of several classical test problems verify high order accuracy of the method. We have also applied the method to a full wave inversion problem. 

We have only used uniform Cartesian grids in the numerical experiments in this paper. For the problems with strong interface phenomena, locally refined mesh would be more efficient. We do not expect difficulties in the extension in this direction, though we have not pursued its implementation in this paper. Other future work includes error analysis to improve the understanding of the suboptimal convergence rates seen in some of our experiments, as well as full wave inversion using the adjoint state method for an improved computational efficiency. 

\section*{Acknowledgement}
Appel\"{o} was supported in part by NSF Grant DMS-1319054. Any conclusions 
or recommendations expressed in this paper are those of the author and do not necessarily reflect the views of the NSF. Wang would like to thank Prof. Gunilla Kreiss at Uppsala University for her  support in this work. 

\bibliography{appelo}%

\begin{thebibliography}{10}

\bibitem{antonietti2018high}
P.~F. Antonietti, F.~Bonaldi, and I.~Mazzieri.
\newblock A high-order discontinuous galerkin approach to the elasto-acoustic
  problem.
\newblock {\em arXiv:1803.01351}, 2018.

\bibitem{Upwind2}
D.~Appel\"o and T.~Hagstrom.
\newblock A new discontinuous {G}alerkin formulation for wave equations in
  second order form.
\newblock {\em {SIAM} Journal on Numerical Analysis}, 53(6):2705--2726, 2015.

\bibitem{el_dg_dath}
D.~Appel{\"{o}} and T.~Hagstrom.
\newblock An energy-based discontinuous {G}alerkin discretization of the
  elastic wave equation in second order form.
\newblock {\em Computer Methods in Applied Mechanics and Engineering},
  338:362--391, 2018.

\bibitem{L-BFGS}
R.~Byrd, P.~Lu, J.~Nocedal, and C.~Zhu.
\newblock A limited memory algorithm for bound constrained optimization.
\newblock {\em SIAM Journal on Scientific Computing}, 16(5):1190--1208, 1995.

\bibitem{ChouShuXing2014}
C.-S. Chou, C.-W. Shu, and Y.~Xing.
\newblock Optimal energy conserving local discontinuous {G}alerkin methods for
  second-order wave equation in heterogeneous media.
\newblock {\em Journal of Computational Physics}, 272:88 -- 107, 2014.

\bibitem{ChungT119}
E.~T. Chung, C.~Y. Lam, and J.~Qian.
\newblock A staggered discontinuous {G}alerkin method for the simulation of
  seismic waves with surface topography.
\newblock {\em Geophysics}, 80(4):T119--T135, 2015.

\bibitem{Duru2014}
K.~Duru and K.~Virta.
\newblock Stable and high order accurate difference methods for the elastic
  wave equation in discontinuous media.
\newblock {\em Journal of Computational Physics}, 279:37 -- 62, 2014.

\bibitem{TR_Givoli}
D.~Givoli.
\newblock Time reversal as a computational tool in acoustics and
  elastodynamics.
\newblock {\em Journal of Computational Acoustics}, 22(03):1430001, 2014.

\bibitem{GSSwave}
M.~J. Grote, A.~Schneebeli, and D.~Sch{\"o}tzau.
\newblock Discontinuous {G}alerkin finite element method for the wave equation.
\newblock {\em SIAM Journal on Numerical Analysis}, 44(6):2408--2431, 2006.

\bibitem{Hagstrom2012}
T.~Hagstrom and G.~Hagstrom.
\newblock Grid stabilization of high-order one-sided differencing {II}:
  {S}econd-order wave equations.
\newblock {\em Journal of Computational Physics}, 231(23):7907 -- 7931, 2012.

\bibitem{KasDum2D06}
M.~K\"{a}ser and M.~Dumbser.
\newblock An arbitrary high-order discontinuous {G}alerkin method for elastic
  waves on unstructured meshes - i. the two-dimensional isotropic case with
  external source terms.
\newblock {\em Geophysical Journal International}, 166(2):855--877, 2006.

\bibitem{Kaser2008}
M.~K\"{a}ser and M.~Dumbser.
\newblock A highly accurate discontinuous {G}alerkin method for complex
  interfaces between solids and moving fluids.
\newblock {\em Geophysics}, 73(3):T23--T35, 2008.

\bibitem{KreOli72}
H.-O. Keiss and J.~Oliger.
\newblock Comparison of accurate methods for the integration of hyperbolic
  equations.
\newblock {\em Tellus}, 24:199--215, 1972.

\bibitem{Komatitsch2000}
D.~Komatitsch, C.~Barnes, and J.~Tromp.
\newblock Wave propagation near a fluid-solid interface: {A} spectral-element
  approach.
\newblock {\em Geophysics}, 65(2):623--631, 2000.

\bibitem{KOP_bc}
H.~O. Kreiss, O.~E. Ortiz, and N.~A. Petersson.
\newblock Initial-boundary value problems for second order systems of partial
  differential equations.
\newblock {\em ESAIM M2AN}, (46):559--593, 2012.

\bibitem{monkola2011numerical}
S.~M{\"o}nk{\"o}l{\"a}.
\newblock Numerical simulation of fluid-structure interaction between acoustic
  and elastic waves.
\newblock {\em PhD Thesis in University of Jyv{\"a}skyl{\"a}}, 2011.

\bibitem{Nguyen2011}
N.~Nguyen, J.~Peraire, and B.~Cockburn.
\newblock High-order implicit hybridizable discontinuous {G}alerkin methods for
  acoustics and elastodynamics.
\newblock {\em Journal of Computational Physics}, 230(10):3695 -- 3718, 2011.

\bibitem{Nocedal2006NO}
J.~Nocedal and S.~J. Wright.
\newblock {\em Numerical Optimization}.
\newblock Springer, New York, 2nd edition, 2006.

\bibitem{Petersson2015}
N.~A. Petersson and B.~Sj\"{o}green.
\newblock Wave propagation in anisotropic elastic materials and curvilinear
  coordinates using a summation-by-parts finite difference method.
\newblock {\em Journal of Computational Physics}, 299:820--841, 2015.

\bibitem{Petersson2018}
N.~A. Petersson and B.~Sj\"{o}green.
\newblock High order accurate finite difference modeling of seismo-acoustic
  wave propagation in a moving atmosphere and a heterogeneous earth model
  coupled across a realistic topography.
\newblock {\em Journal of Scientific Computing}, 74(1):290--323, 2018.

\bibitem{DORMAND198019}
P.~Prince and J.~Dormand.
\newblock High order embedded {R}unge-{K}utta formulae.
\newblock {\em Journal of Computational and Applied Mathematics}, 7(1):67 --
  75, 1981.

\bibitem{Sheldon2016}
J.~P. Sheldon, S.~T. Miller, and J.~S. Pitt.
\newblock A hybridizable discontinuous {G}alerkin method for modeling
  fluid-structure interaction.
\newblock {\em Journal of Computational Physics}, 326:91--114, 2016.

\bibitem{Stanglmeier2016}
M.~Stanglmeier, N.~Nguyen, J.~Peraire, and B.~Cockburn.
\newblock An explicit hybridizable discontinuous {G}alerkin method for the
  acoustic wave equation.
\newblock {\em Computer Methods in Applied Mechanics and Engineering}, 300:748
  -- 769, 2016.

\bibitem{Virta2014}
K.~Virta and K.~Mattsson.
\newblock Acoustic wave propagation in complicated geometries and heterogeneous
  media.
\newblock {\em Journal of Scientific Computing}, 61(1):90--118, 2014.

\bibitem{Wang2016}
S.~Wang, K.~Virta, and G.~Kreiss.
\newblock High order finite difference methods for the wave equation with
  non-conforming grid interfaces.
\newblock {\em Journal of Scientific Computing}, 68(3):1002--1028, 2016.

\bibitem{Wilcox2010}
L.~C. Wilcox, G.~Stadler, C.~Burstedde, and O.~Ghattas.
\newblock A high-order discontinuous {G}alerkin method for wave propagation
  through coupled elastic-acoustic media.
\newblock {\em Journal of Computational Physics}, 229(24):9373 -- 9396, 2010.

\bibitem{ruichao}
R.~Ye, M.~V. de Hoop, C.~Petrovitch, L.~Pyrak-Nolte, and L.~C. Wilcox.
\newblock A discontinuous galerkin method with a modified penalty flux for the
  propagation and scattering of acousto-elastic waves.
\newblock {\em Geophysical Journal International}, 205(2):1267--1289, 2016.

\end{thebibliography}
\bibliographystyle{abbrv}

\end{document}